\documentclass[11pt, a4paper]{article}

\usepackage[english]{babel}
\usepackage{dutchcal}

\usepackage{amsmath}
\usepackage{amssymb}
\usepackage{amsthm}
\usepackage{tikz}
\usepackage{tikz-cd}
\usepackage{todonotes}
\usepackage{hyperref}
\usepackage{cleveref}
\usepackage{footmisc}
\usepackage{mathrsfs}
\usepackage{extarrows}
\usepackage{enumitem}
\usepackage{stmaryrd}

\newcommand{\llpar}{\llparenthesis}
\newcommand{\rrpar}{\rrparenthesis}

\newcommand{\ab}{\mathrm{ab}}

\newcommand{\afp}{\mathrm{afp}}

\newcommand{\Ell}{\mathrm{Ell}}
\newcommand{\QEll}{\mathrm{QEll}}
\newcommand{\gQEll}{\mathbf{QEll}}

\newcommand{\Temp}{\mathrm{Temp}}

\newcommand{\fin}{\mathrm{fin}}

\newcommand{\gl}{\mathrm{gl}}
\newcommand{\km}{\mathrm{KM}}
\newcommand{\GHM}{\mathrm{GHM}}

\newcommand{\ori}{\mathrm{or}}

\newcommand{\Tate}{\mathrm{Tate}}
\renewcommand{\top}{\mathrm{top}}

\newcommand{\temp}{\mathrm{temp}}
\newcommand{\fingl}{\fin\textup{-}\gl}

\DeclareMathOperator{\MF}{MF}

\DeclareMathOperator{\KU}{KU}
\DeclareMathOperator{\ku}{ku}
\newcommand{\gKU}{\mathbf{KU}}
\DeclareMathOperator{\KO}{KO}
\newcommand{\gKO}{\mathbf{KO}}

\newcommand{\Sph}{\mathbf{S}}
\newcommand{\gTMF}{\mathbf{TMF}}

\DeclareMathOperator{\TMF}{TMF}
\DeclareMathOperator{\tmf}{tmf}
\DeclareMathOperator{\Tmf}{Tmf}

\DeclareMathOperator{\Ab}{Ab}

\DeclareMathOperator{\AVar}{AVar}
\newcommand{\BTg}{{\mathrm{BT}_{2g}}}
\newcommand{\BTn}{{\mathrm{BT}_n}}
\newcommand{\BTtwo}{{\mathrm{BT}_2}}
\newcommand{\BTtwoor}{{\mathrm{BT}_2^\ori}}
\DeclareMathOperator{\CAlg}{CAlg}

\DeclareMathOperator{\Cyc}{Cyc}
\DeclareMathOperator{\CycEll}{CycEll}

\DeclareMathOperator{\Fun}{Fun}

\DeclareMathOperator{\Glo}{Glo}

\DeclareMathOperator{\Mod}{Mod}

\DeclareMathOperator{\Lat}{Lat}

\newcommand{\Ori}{\mathrm{Or}}
\DeclareMathOperator{\QCoh}{QCoh}

\DeclareMathOperator{\Spec}{Spec}

\DeclareMathOperator{\Shv}{Shv}
\DeclareMathOperator{\Stk}{Stk}
\DeclareMathOperator{\Sp}{Sp}

\newcommand{\Spc}{\mathcal{S}}

\newcommand{\colim}{\mathrm{colim}\,}

\newcommand{\id}{\mathrm{id}}
\DeclareMathOperator{\Map}{Map}

\DeclareMathOperator{\Hom}{Hom}
\newcommand{\HOM}{\underline{\mathrm{Hom}}}
\newcommand{\op}{\mathrm{op}}

\newcommand{\A}{\mathsf{A}}
\newcommand{\B}{\mathbf{B}}
\renewcommand{\b}{\mathsf{B}}
\newcommand{\BT}{\mathbf{BT}}

\newcommand{\C}{\mathcal{C}}

\newcommand{\E}{\mathbf{E}}
\newcommand{\e}{\mathsf{E}}
\newcommand{\EE}{\mathscr{E}}

\newcommand{\G}{\mathbf{G}}

\newcommand{\M}{\mathsf{M}}
\newcommand{\N}{\mathbf{N}}

\renewcommand{\O}{\mathcal{O}}

\renewcommand{\P}{\mathbf{P}}
\newcommand{\Q}{\mathbf{Q}}

\newcommand{\QQ}{\mathcal{Q}}
\newcommand{\R}{\mathbf{R}}

\renewcommand{\SS}{\mathfrak{S}}
\newcommand{\s}{\mathsf{X}}
\renewcommand{\t}{\mathsf{T}}
\newcommand{\T}{\mathrm{T}}

\newcommand{\TT}{\mathbf{T}}

\newcommand{\gX}{\mathbf{X}}

\newcommand{\x}{\mathsf{X}}

\newcommand{\y}{\mathsf{Y}}

\newcommand{\Z}{\mathbf{Z}}

\newcommand{\ga}{\gamma}
\newcommand{\Ga}{\Gamma}

\usepackage{mathtools}

\theoremstyle{theorem}\numberwithin{equation}{section}
\newtheorem{theorem}[equation]{Theorem}
\crefname{theorem}{{th}.\!\!}{{ths}.\!\!}
\Crefname{theorem}{{Th}.\!\!}{{Ths}.\!\!}
\newtheorem{theoremalph}{Theorem}

\newtheorem{coralph}[theoremalph]{Corollary}

\crefname{coralph}{{cor}.\!\!}{{cors}.\!\!}
\Crefname{coralph}{{Cor}.\!\!}{{Cors}.\!\!}
\crefname{theoremalph}{{th}.\!\!}{{ths}.\!\!}
\Crefname{theoremalph}{{Th}.\!\!}{{Ths}.\!\!}

\Crefname{problem}{{Prb}.\!\!}{{Prbs}.\!\!}
\newtheorem{prop}[equation]{Proposition}
\Crefname{prop}{{Pr}.\!\!}{{Prs}.\!\!}
\newtheorem{lemma}[equation]{Lemma}
\Crefname{lemma}{{Lm}.\!\!}{{Lms}.\!\!}
\newtheorem{cor}[equation]{Corollary}
\Crefname{cor}{{Cor}.\!\!}{{Cors}.\!\!}

\Crefname{conjecture}{{Conj}.\!\!}{{Conjs}.\!\!}

\theoremstyle{definition}\numberwithin{equation}{section}
\newtheorem{mydef}[equation]{Definition}
\Crefname{mydef}{{Df}.\!\!}{{Dfs}.\!\!}

\Crefname{recall}{{Rcl}.\!\!}{{Rcls}.\!\!}
\newtheorem{construction}[equation]{Construction}
\Crefname{construction}{{Con}.\!\!}{{Cons}.\!\!}

\Crefname{ass}{{As}.\!\!}{{As}.\!\!}

\Crefname{notation}{{Nt}.\!\!}{{Nts}.\!\!}

\Crefname{situation}{{St}.\!\!}{{Sts}.\!\!}

\theoremstyle{remark}\numberwithin{equation}{section}
\newtheorem{example}[equation]{Example}
\Crefname{example}{{Ex}.\!\!}{{Exs}.\!\!}

\Crefname{nonexample}{{NonEx}.\!\!}{{NonEx}.\!\!}

\Crefname{claim}{{Clm}.\!\!}{{Clms}.\!\!}
\newtheorem{remark}[equation]{Remark}
\Crefname{remark}{{Rmk}.\!\!}{{Rmks}.\!\!}

\Crefname{idea}{{Id}.\!\!}{{Ids}.\!\!}
\newtheorem{warn}[equation]{Warning}
\Crefname{warn}{{Warn}.\!\!}{{Warns}.\!\!}

\Crefname{question}{{Qn}.\!\!}{{Qns}.\!\!}

\Crefname{figure}{{Fig.}\!\!}{{Figs.}\!\!}
\Crefname{footnote}{{Fn.}\!\!}{{Fn.}\!\!}

\Crefname{part}{{\textsection}\!\!}{{\textsection}\!\!}
\Crefname{chapter}{{\textsection}\!\!}{{\textsection}\!\!}
\Crefname{section}{{\textsection}\!\!}{{\textsection}\!\!}
\Crefname{subsection}{{\textsection}\!\!}{{\textsection}\!\!}
\Crefname{appendix}{{\textsection}\!\!}{{\textsection}\!\!}

\mathtoolsset{centercolon}

\newlist{numberenum}{enumerate}{1}
\setlist[numberenum]{label={{\upshape(\arabic*)}}}

\begin{document}
\title{On the derived Tate curve and global smooth Tate $K$-theory}

\author{Jack Morgan Davies\footnote{\href{mailto:davies@uni-wuppertal.de}{\texttt{davies@uni-wuppertal.de}}}\, and Sil Linskens\footnote{\href{mailto:sil.linskens@mathematik.uni-regensburg.de}{\texttt{sil.linskens@mathematik.uni-regensburg.de}}}}
\maketitle

\begin{abstract}
    The interplay between equivariant stable homotopy theory and spectral algebraic geometry is used to construct a derived Tate curve over $\KU\llpar q\rrpar$, a lift of the classical elliptic curve of Tate over $\Z\llpar q \rrpar$. Applications of both an algebro-geometric and a topological flavour follow. First, we construct a spectral algebro-geometric model for the compactification of the moduli stack of oriented elliptic curves, giving a canonical choice of holomorphic topological $q$-expansion map. Then we define globally equivariant forms of Tate $K$-theory $\gKO\llpar q \rrpar$ and $\gKU\llpar q \rrpar$, and equip them with globally equivariant meromorphic topological $q$-expansion maps from global topological modular forms. Finally, we explore $C_2$-equivariant versions of global Tate $K$-theory and connect them with $C_2$-equivariant global topological modular forms with level structures.
\end{abstract}

\setcounter{tocdepth}{2}
\tableofcontents

\newpage

\section{Introduction}
The classical Tate curve $\T^\heartsuit$ over $\Z\llpar q \rrpar$ plays a central role in the study of elliptic curves and modular forms. For example, the map $\Spec \Z\llpar q \rrpar \to \M_\Ell$ into the moduli stack of elliptic curves defined by $\T^\heartsuit$ induces the $q$-expansion map on global sections
\begin{equation}\label{qexpmap}
    \MF \to \Z\llpar q \rrpar, \qquad f \mapsto f(\tau) = \sum_{n}a_n(f) q^n,
\end{equation}
sending a meromorphic modular form to its $q$-expansion.\\

The goal of this article is to enhance the Tate curve from an elliptic curve over $\pi_0 \KU\llpar q\rrpar = \Z\llpar q \rrpar$ to a derived elliptic curve $\T$ over the $\E_\infty$-ring $\KU\llpar q \rrpar$ (\Cref{existenceoftatecurve}). We then use this derived elliptic curve to compactify the derived moduli stack of oriented elliptic curves (\Cref{comparison}) and produce a highly structured refinement of (\ref{qexpmap}) to a morphism $\gTMF \to \gKO\llpar q \rrpar$ of \emph{global $\E_\infty$-rings} between equivariant topological modular form and equivariant Tate $K$-theory (\Cref{globalmap}). Variations on this theme are also applied to define \emph{global quasi-elliptic cohomology} and forms of \emph{$C_2$-equivariant Tate $K$-theory} and prove some of their basic properties (\Cref{cor:qellstuffintro,cor:realstuffintro}). The techniques used to construct the derived Tate curve $\T$ combine concepts from equivariant stable homotopy theory and spectral algebraic geometry. In particular, it involves the \emph{cyclification} of equivariant cohomology theories, an application of Lurie's Serre--Tate theorem for derived elliptic curves, and a corepresentability property of $\TMF$. As a whole, this serves to further strengthen the relationship between arithmetic geometry, spectral algebraic geometry, and equivariant homotopy theory.

\subsection*{Background}

The universal elliptic cohomology theory $\TMF$ of \emph{topological modular forms} was first constructed by Hopkins and collaborators \cite{hopkinsfirsttmficm}. Since its inception, this cohomology theory has been a bridge between homotopy theory and other areas such as number theory, differential geometry, and physics (\cite{hopkinssecondtmficm,stolzteichner}). The spectrum $\TMF$ has also played a transformative role in stable homotopy theory itself (\cite{lurieecsurvey}).\\

The original construction of $\TMF$ proceeded via obstruction theory; see \cite[\textsection12]{tmfbook}, for example. More recently, Lurie has given an elegant construction of $\TMF$ using spectral algebraic geometry. To do this, Lurie first defines a notion of an \emph{oriented elliptic curves} over an $\E_\infty$-ring and in general over derived stacks. One then proves that the moduli problem for such elliptic curves is represented by the derived analog of a Deligne--Mumford stack $\M_\Ell^\ori$, and defines the $\E_\infty$-ring $\TMF$ as the global sections of this stack. Far from just being another construction of $\TMF$, Lurie's approach has many other benefits. For example, this construction yields a \emph{universal oriented elliptic curve} $\e$ over $\M_\Ell^\ori$, a key ingredient in the construction of \emph{equivariant topological modular forms} $\gTMF$ (\cite{lurieecsurvey, davidandlennart,ec3}), itself a derived refinement of Grojnowski's influential construction of equivariant elliptic cohomology over the complex numbers (\cite{Grojnowskielliptic}).\\

As this demonstrates, there is much to be gained from constructing elliptic cohomology theories together with derived lifts of classical elliptic curves. The Tate curve over $\Z\llpar q \rrpar$ is a particularly important candidate to find a derived refinement of, given its importance to number theory and arithmetic geometry, physics and the Witten genus, and as a valuable connection between topological modular forms and topological $K$-theory. Despite recent advances though, it remains a difficult task to construct oriented elliptic curves. For example, no technology is available to define elliptic curves over $\E_\infty$-rings, such as the Tate curve, via Weierstraß equations. Even if one could, the group structure on an oriented elliptic curve is \emph{extra data} and does not simply follow from the \emph{property} that a derived scheme is a smooth curve of genus one with a chosen point. As another approach, one might try to use Lurie's \emph{étale rigidity theorem}, which in this case states that the $\infty$-category of derived stacks which are étale over $\M_\Ell^\ori$ is equivalent to the $2$-category of classical Deligne--Mumford stacks which are étale over the classical moduli stack of elliptic curves $\M_\Ell$. This produces canonical oriented elliptic curves over moduli stacks of elliptic curves with level structure $\M^\ori_1(n)$ and $\M^\ori_0(n)$, for example. However, no map $\Spec \Z\llpar q \rrpar \to \M_\Ell$ could be étale as it cannot be of relative dimension zero. As a result, we see no methods to purely formally construct a derived Tate curve.

\subsection*{Main result}
In this article, we produce a lift of the Tate curve to an oriented elliptic curve over the $\E_\infty$-ring $\KU\llpar q \rrpar$. The properties of this refinement are summarized in the following theorem, a slightly refined version of which appears as \Cref{liftoftate}.

\begin{theoremalph}\label{existenceoftatecurve}
There exists a unique oriented elliptic curve $\T$ over $\KU\llpar q\rrpar$ such that:
\begin{enumerate}
    \item its underlying elliptic curve is the classical Tate curve $\T^\heartsuit$ over $\Z\llpar q\rrpar$,
    \item its associated formal group is the formal multiplicative group scheme $\widehat{\G}_{m}$ over $\KU\llpar q \rrpar$,
    \item the rationalisation of the associated map of $\E_\infty$-rings $\TMF \to \KU\llpar q \rrpar$, a \emph{meromorphic topological $q$-expansion map}, is the rationalisation of (\ref{qexpmap}), and
    \item its associated oriented $\P$-divisible group is the base-change of $\T_\km$ of \Cref{spectralliftofkm} from $\KU[q^\pm]$ to $\KU\llpar q\rrpar$, where $\T_\km$ is a spectral lift of the Katz--Mazur group scheme of \cite[\textsection 8.7]{km} defined over $\Z[q^\pm]$.
\end{enumerate}
\end{theoremalph}

In a complementary fashion to how spectral algebro-geometric techniques are crucial to constructing equivariant topological modular forms, our construction of the oriented elliptic curve $\T$ relies heavily on equivariant homotopy theory. Essentially, we first construct the equivariant cohomology theory (for finite groups) associated with $\T$ by \emph{cyclifying} equivariant complex topological $K$-theory (\Cref{emergenceofpdivgroup}). This immediately defines a preoriented $\P$-divisible group $\T_\km$ which is a spectral lift of the Katz--Mazur group scheme from $\Z[q^\pm]$ to $\KU[q^\pm]$ (\Cref{spectralliftofkm}). Then spectral algebro-geometric technology kicks into play. Using a variant of Lurie's spectral Serre--Tate theorem (\Cref{spectralserretate}), we can combine $\T_\km$ over $\KU\llpar q \rrpar^\wedge_p$ together with the classical Tate curve over $\Z\llpar q \rrpar^\wedge_p$ to obtain an elliptic curve $\widehat{\T}$ over the $\KU\llpar q \rrpar^\wedge_p$ for each prime $p$. These elliptic curves are then glued together with some rational information using the corepresentability of $\TMF$ into the global sections of well-behaved stacks (\Cref{corepresentabilityoftmfinnicecases}) to obtain an elliptic curve over $\KU\llpar q \rrpar$.

\subsection*{Applications}

Once we have constructed the oriented Tate curve $\T$, we apply it in two orthogonal directions: firstly to compactify the spectral moduli stack of oriented elliptic curves (\Cref{comparison}), and secondly to explore refinements of the meromorphic $q$-expansion map further to a map $\gTMF \to \gKO\llpar q \rrpar$ of \emph{global $\E_\infty$-rings} (\Cref{globalmap}) as well as some other equivariant variations (\Cref{cor:qellstuffintro,cor:realstuffintro}).

\subsubsection*{Compactifying the moduli stack of oriented elliptic curves}

Noting that the automorphism group of $\T$ is generated by the inversion map $[-1]$, we show that $\T$ descends to an oriented elliptic curve over $\Spec \KU\llpar q \rrpar / C_2$, where $C_2$ acts on $\KU\llpar q \rrpar$ by complex-conjugation. In particular, the global sections of $\Spec \KU\llpar q \rrpar /C_2$ can be identified with $\KO\llpar q \rrpar$. After taking global sections and homotopy groups, the map $\Spec \KU\llpar q \rrpar/C_2 \to \M_\Ell^\ori$ is rationally equivalent to (\ref{qexpmap}). Just as in the classical case, we can now define the \emph{compactification of the moduli stack of oriented elliptic curves}, denoted as $\overline{\M}_\Ell^\ori$, by the pushout diagram of (1-localic) derived stacks
\[\begin{tikzcd}
    {\Spec \KU\llpar q \rrpar /C_2}\ar[d,"{\T}"]\ar[r]    &   {\Spec \KU\llbracket q\rrbracket/C_2}\ar[d] \\
    {\M_{\Ell}^\ori}\ar[r]  &   {\overline{\M}_{\Ell}^\ori.}
\end{tikzcd}\]

\begin{coralph}\label{comparison}
The derived stack $\overline{\M}_\Ell^\ori$ is proper over $\Spec \Sph$ and has underlying classical stack the Deligne--Mumford compactification of the moduli stack of elliptic curves. Moreover, its structure sheaf defines elliptic cohomology theories when restricted to affines. In particular, there is an equivalence of derived stacks
\[\overline{\M}_{\Ell}^\ori\simeq (\overline{\M}_\Ell,\O^\top)\]
where $\O^\top$ is the Goerss--Hopkins--Miller sheaf constructed by Behrens in \cite[\textsection12]{tmfbook}.
\end{coralph}

As an immediate application of this equivalence of derived stacks, we obtain an equivalence of $\E_\infty$-rings $\Ga(\overline{\M}_\Ell^\ori)\simeq \Tmf$. In fact, as the construction of $\overline{\M}_\Ell^\ori$ is done purely in spectral algebraic geometry, we simply define $\Tmf$ as the global sections $\Ga(\overline{\M}_\Ell^\ori)$. From this definition, we obtain a canonical \emph{holomorphic topological $q$-expansion map}
\begin{equation}\label{holoqexp}\Tmf\to \KO \llbracket q \rrbracket,\end{equation}
which induces the holomorphic analogue of (\ref{qexpmap}) on rational homotopy groups. In the literature, this holomorphic topological $q$-expansion map is constructed using obstruction theory. To the best of our knowledge, our construction of this map in \Cref{comparison} is the first to use purely spectral algebraic geometry. In \cite[\textsection A]{hilllawson}, it is shown that (\ref{holoqexp}) is unique up to homotopy, so our definition above agrees with others found in the literature. The formal properties of $\overline{\M}^\ori_\Ell$ listed above also allow one to compute the homotopy groups of these global sections using a descent spectral sequence, as done in \cite{smfcomputation} for example, or determine abstract properties about the stack $\overline{\M}_\Ell^\ori$ such as $0$-affineness, see \cite{akhilandlennart} or \cite{reconstruction}.

\subsubsection*{An equivariant meromorphic topological \texorpdfstring{$q$}{q}-expansion map for compact Lie groups}

Not only does the oriented Tate curve $\T$ give the meromorphic topological $q$-expansion map $\TMF \to \KO\llpar q \rrpar$ as a morphism of $\E_\infty$-rings, but using recent advances in equivariant homotopy theory, this can be refined further to a morphism of \emph{global $\E_\infty$-rings}.\\

In \cite{davidandlennart}, Gepner--Meier construct equivariant cohomology theories associated with oriented elliptic curves, such as the Tate curve $\T$ appearing in \Cref{existenceoftatecurve}. Roughly speaking, for an oriented elliptic curve $\e$ over a derived stack $\s$, they define $\Ell_{\e}(\B G)$ to be the stack of homomorphisms $\HOM(\underline{G}^\vee, \e)$ for each abelian compact Lie group $G$. The value of $\Ell_{\e}(-)$ on a general $G$-space $Y$ for a general compact Lie group $G$ is obtained by appropriately Kan extending this diagram.\\

In particular, this uniformly produces a $G$-equivariant cohomology theory for all compact Lie groups $G$ recovering what is known as  \emph{global spectrum}. Such an object is roughly the data of a $G$-equivariant spectrum $\gX_G$ for each compact Lie group $G$, together with maps $\varphi^\ast(\gX_K)\to \gX_G$ for each map of compact Lie groups $\varphi\colon G\to K$, which are further equivalences if $\varphi$ is injective, plus higher coherence data. This intuition is made precise in \cite{denissilluca}, where it is shown that Schwede's $\infty$-category of global spectra $\Sp_\gl$ can be written as a partially lax limit of $\infty$-categories of $G$-spectra $\Sp_G$. The recent work of Gepner, Pol and the second-named author \cite{gepner2024global2ringsgenuinerefinements} shows that these equivariant cohomology theories $\Ell_{\e}{(-)}$ refine to a \emph{global $\E_\infty$-ring} $\Ga_\gl \Ell_{\e}$, so to an $\E_\infty$-object in $\Sp^\gl$. One can apply $\Ga_\gl \Ell_{(-)}$ to the universal oriented elliptic curve $\e^\ori$ living of the moduli stack of oriented elliptic curves $\M_\Ell^\ori$ to obtain \emph{global topological modular forms} $\gTMF$; see \cite[Ex.14.22]{gepner2024global2ringsgenuinerefinements}. Similarly, feeding the derived Tate curve $\T$ over either $\Spec \KU\llpar q \rrpar$ or $\Spec \KU\llpar q \rrpar / C_2$ into $\Ga_\gl \Ell_{(-)}$ leads to the global $\E_\infty$-rings $\gKU\llpar q \rrpar$ and $\gKO\llpar q \rrpar$ of \emph{global smooth real} and \emph{complex Tate $K$-theory}, respectively. The adjective \emph{smooth} refers to the fact that the Tate curve over $\Spec \KU\llpar q \rrpar / C_2$ is a smooth elliptic curve, rather than a potential curve with nodal singularity over $\Spec \KU\llbracket q \rrbracket$ that we do not consider here. One advantage of this global refinement is that these global $\E_\infty$-rings come with a coherent collection of transfer maps, which is a piece of structure not immediately visible from the stacks $\Ell_{\e}(\B G)$.\\

One of the inspirations of this article was to construct a globally equivariant version of the meromorphic topological $q$-expansion map appearing in part 3 of \Cref{existenceoftatecurve}. This arises by combining our construction of $\T$ with the aforementioned works.

\begin{coralph}\label{globalmap}
There exists a \emph{global meromorphic topological $q$-expansion map}
\[\gTMF\to \gKO\llpar q\rrpar\] 
in the $\infty$-category $\CAlg(\Sp^\gl)$, whose underlying map of $\E_\infty$-rings is that of part 3 in \Cref{existenceoftatecurve}.
\end{coralph}

In particular, for any global space $X$, we have a map of rings
\[\gTMF^\ast(X)\to \gKO\llpar q \rrpar^\ast(X).\]
For instance, if $Y$ is a topological space with a continuous $G$-action by a compact Lie group $G$, then we can set $X = Y // G$ and the above yields a natural multiplicative map between $G$-equivariant cohomology theories
\[\gTMF_G^\ast(Y) \to \gKO\llpar q\rrpar_G^\ast(Y).\]
As these maps of rings come from a map of global $\E_\infty$-rings, they also preserve (naïve) power operations, such as Dyer--Lashoff operations and the Tate-valued Frobenius, as well as commute with all transfer maps along monomorphisms of groups.\\

The Gepner--Meier approach to equivariant elliptic cohomology works particularly well when applied to compact Lie groups. If one is only interested in \emph{finite} groups, ie, in fin-global spectra, then one should rather work in the context of Lurie's \emph{tempered cohomology theories}. The question of if such theories are represented by global spectra is also solved in \cite{gepner2024global2ringsgenuinerefinements}, however, we prefer to use the more general and flexible approach of the authors and Balderrama in \cite{tempered}. These equivariant theories only depend on the underlying oriented $\P$-divisible group of an oriented elliptic curve, so we may as well work with the spectral lift of the Katz--Mazur $\P$-divisible group $\T_\km$ over the more initial $\E_\infty$-ring $\KU[q^\pm]$.

\subsubsection*{Global quasi-elliptic cohomology}
The fin-global $\E_\infty$-ring associated with $\T_\km$ is denoted by $\gQEll$ and called \emph{global quasi-elliptic cohomology}. Elsewhere, such as \cite[Ex.5.1.6]{tempered}, we denote $\gQEll$ by $\gKU_\km$, a kind of \emph{Katz--Mazur $K$-theory}. The following is a summary and simplification of \Cref{basechangeprop} and \Cref{kuenneth}.

\begin{coralph}\label{cor:qellstuffintro}
    There is an equivalence of fin-global $\E_\infty$-rings
    \[\gQEll \otimes_{\KU[q^\pm]} \KU\llpar q \rrpar \simeq \gKU_\fin\llpar q \rrpar,\]
    where $\gKU_\fin\llpar q \rrpar$ is the restriction of global smooth Tate $K$-theory to a fin-global $\E_\infty$-ring. For finite groups $H,G$, the natural map of graded rings
    \[\gQEll^\ast_H \underset{\Z[q^\pm, u^\pm]}{\otimes} \gQEll^\ast_G \xrightarrow{\simeq} \gQEll^\ast_{H\times G}\]
    is an isomorphism.
\end{coralph}

By \Cref{rmk:compwithgrh}, the coefficients of $\gQEll$ agree with the equivariant quasi-elliptic cohomology of Ganter--Rezk--Huan. These properties of $\gQEll$ are coherent versions of similar isomorphisms found in Huan's work \cite{huanquasielliptic}, and their proof outlines how general theorems of tempered cohomology theories apply to $\gQEll$.

\subsubsection*{Two \texorpdfstring{$C_2$}{C2}-equivariant global meromorphic topological \texorpdfstring{$q$}{q}-expansion maps}
Our final application of the derived Tate curve $\T$, is to show how it can be used to construct further analogues of the global meromorphism topological $q$-expansion map of \Cref{globalmap}.\\

There is a natural $C_2$-action on $\KU\llpar q \rrpar$ coming from complex conjugation and there are also natural $C_2$-actions on various forms of topological modular forms with level structure such as $\TMF_1(3)$. Our final goal is to show how a further understanding of $\T$ and $\T_\km$ can lead to $C_2$-equivariant refinements of the global meromorphic $q$-expansion maps of \Cref{globalmap}. To this end, we write $\gKU\llpar q \rrpar (\omega)$ for the $\gKU\llpar q \rrpar$-algebra given by (implicitly) inverting $3$ and adjoining a primitive $3$\textsuperscript{rd} root of unity, and $\gKU\llpar q^{1/3} \rrpar[\tfrac{1}{3}]$ for the globally flat extension of $\gKU\llpar q \rrpar$ which adjoins a $3$\textsuperscript{rd} root of $q$ and inverts $3$. The appropriate $C_2$-actions are defined in \Cref{ssec:real}. The following is a summary of \Cref{realunramqexpansionmap,realramfqexpansionmap}.

\begin{coralph}\label{cor:realstuffintro}
    There are maps in $\Fun(BC_2, \CAlg(\Sp_\gl))$ of the form
    \[\ga\colon \gTMF_1(3) \to \gKU\llpar q \rrpar(\omega), \qquad \rho \colon \gTMF_1(3) \to \gKU\llpar q^{1/3}\rrpar[\tfrac{1}{3}]\]
    such that the underlying diagram of global $\E_\infty$-rings
    \[\begin{tikzcd}
        {\gTMF_1(3)}\ar[d, "\ga"]   &   {\gTMF_1(3)^{hC_2}}\ar[l]\ar[d, "{\ga^{hC_2}}"]    &   {\gTMF}\ar[l]\ar[d, "\T"]\ar[r]  &   {\gTMF_1(3)^{hC_2}}\ar[r]\ar[d, "{\rho^{hC_2}}"]    &   {\gTMF_1(3)}\ar[d, "\rho"]   \\
        {\gKU\llpar q \rrpar(\omega)}   &   {\gKU\llpar q \rrpar(\omega)^{hC_2}}\ar[l]  &   {\gKO\llpar q \rrpar}\ar[r]\ar[l]    &   {\gKU\llpar q^{1/3} \rrpar[\tfrac{1}{3}]^{hC_2}}\ar[r]  &   {\gKU\llpar q^{1/3} \rrpar[\tfrac{1}{3}]}
    \end{tikzcd}\]
    commutes, where $\T$ is the map of \Cref{globalmap} induced by the derived Tate curve $\T$ over $\Spec \KU\llpar q\rrpar/C_2$.
\end{coralph}

Taking underlying $\E_\infty$-rings, one obtains maps between Borel $C_2$-equivariant $\E_\infty$-rings
    \[\KU\llpar q \rrpar(\omega) \gets \TMF_1(3) \to \KU\llpar q^{1/3}\rrpar[\tfrac{1}{3}],\]
the middle of which was extensively analysed in \cite{lennartandhill}. This yet again highlights how having the derived elliptic curve $\T$ can be used to obtain more structured connections between topological modular forms and topological $K$-theory.


\subsection*{Outline}

As mentioned, we first construct the equivariant cohomology theory associated with the derived Tate curve. To do so, we set up a procedure for producing interesting new equivariant cohomology theories from old ones. This comes in the form of the \emph{cyclification} endofunctor $\Cyc$ on global spaces in \Cref{sec:cyc}. In \Cref{sec:cyclifyell}, we show that precomposing an equivariant elliptic cohomology theory with the cyclification functor again yields an equivariant cohomology theory. In \Cref{kmsection}, we show that precomposing $\gKU$ with the cyclification functor yields a spectral lift $\T_\km$ of the Katz--Mazur group scheme from $\Z[q^\pm]$ to $\KU[q^\pm]$. In \Cref{tatecurvesection}, we combine $\T_\km$ with a variant of Lurie's spectral Serre--Tate theorem and a corepresentability fact about $\TMF$ from the author's work \cite{reconstruction} with Balderrama to obtain the desired derived Tate curve $\T$ and prove \Cref{existenceoftatecurve}.\\

We then give two broad applications of \Cref{existenceoftatecurve}. In \Cref{sec:compact}, we use $\T$ to compactify the spectral moduli stack of oriented elliptic curves $\M_\Ell^\ori$ and prove \Cref{comparison}. In \Cref{sec:globalthings}, we prove \Cref{globalmap,cor:qellstuffintro,cor:realstuffintro} in three subsections. First, in \Cref{ssec:globaltate}, we use $\T$ and the universal oriented elliptic curve $\e^\ori$ over $\M_\Ell^\ori$ to define a map of global $\E_\infty$-rings $\gTMF \to \gKO\llpar q \rrpar$ and prove \Cref{globalmap}. Then in \Cref{ssec:qell}, we use our spectral lift of the Katz--Mazur group scheme $\T_\km$ over $\KU[q^\pm]$ to define a global version quasi-elliptic cohomology and prove \Cref{cor:qellstuffintro}. Lastly, in \Cref{ssec:real}, we initiate the study of $C_2$-equivariant topological $q$-expansion maps, connecting the $C_2$-equivariant forms of topological modular forms found in \cite{lennartandhill} with forms of Atiyah's Real $K$-theory of \cite{ktheoryandreality}. In particular, this proves \Cref{cor:realstuffintro}.


\subsection*{Previous and future work}
The combination of the Tate curve, equivariant Tate $K$-theory, and something similar to the cyclification functor $\Cyc$ has made many appearances in the literature; see \cite{ganterstringy,rezkquasielliptic,huanquasielliptic}. One highlight of these works is the ability to study equivariant power operations, a topic we would like to come back to in future work. We have also highlighted the construction of a geometric object $\T$ over $\KU\llpar q \rrpar$ and how that leads to a canonical definition of $\overline{\M}_\Ell^\ori$ and $\Tmf$, however, we do not produce a geometric object over $\overline{\M}_\Ell^\ori$. Essentially, this is because we are working with abelian group objects such as elliptic curves, and the classical objects classified by $\overline{\M}_\Ell$ are \emph{generalised elliptic curves}, which do not generally have a group structure. Nevertheless, Lurie suggests a method to construct a geometric object over $\overline{\M}_\Ell^\ori$ by constructing a Tate curve over $\Sph\llbracket q \rrbracket$ using toric spectral algebraic geometry; Laurent Smits has communicated to us that he is carrying out this project in detail. The uniqueness statement of \Cref{existenceoftatecurve} seems useful to compare our two \emph{a priori} different constructions of $\T$ over $\KU\llpar q \rrpar$.


\subsection*{Notation}
The language of higher categories, higher algebra, and spectral algebraic geometry will be used freely; this is summarised in \cite[p.9-11]{ec2}, for example. In particular, we will call $\infty$-categories simply \emph{categories}, and qcqs nonconnective spectral Deligne--Mumford stacks will be called \emph{stacks} and the associated category written as $\Stk$ (in contrast to much of the literature); the word \emph{derived} was only used in the introduction as a catch-all for such stacks. The subcategory of such stacks $\x$ which are \emph{$1$-stacks}, meaning that the space $\Map_{\Stk}(\Spec R, \x)$ is $1$-truncated for each discrete ring $R$, will be called \emph{classical stacks} following \cite[Rmk.1.4.8.3 \& Lm.1.6.8.8]{sag}. The moduli stack of \emph{oriented elliptic curves} will be written as $\M_\Ell^\ori$, see \cite[\textsection 7]{ec2}, which has underlying classical stack $\M_\Ell$, see \cite{delignemumford}. As the only elliptic curves over $\E_\infty$-ring or stacks will be \emph{strict}, see \cite[\textsection1.5]{ec1}, we will simply call them \emph{elliptic curves}.\\ 

We will also use Gepner--Meier's construction of equivariant elliptic cohomology theories from \cite{davidandlennart} frequently; by the main result of \cite{elltempcomp}, this agrees with Lurie's construction of tempered cohomology \cite{ec3} when restricted to $\Spc_{\ab\fin}$. Our notation for global spaces and the global orbit category follow \cite{denissilluca} and are reviewed in \Cref{sec:cyc}.


\subsection*{Acknowledgements}

Thank you to William Balderrama, Anton Engelmann, David Gepner, Lennart Meier, Fabio Neugebauer, Lucas Piessevaux, and Laurent Smits for enlightening conversations on this topic and to Tyler Lawson for a very helpful email exchange.\\

The first author is an associate member of the Hausdorff Center for Mathematics at the University of Bonn (\texttt{DFG GZ 2047/1}, project ID \texttt{390685813}). He would also like to thank the Isaac Newton Institute for Mathematical Sciences, Cambridge, for support and hospitality during the programme \emph{Equivariant homotopy theory in context} where work on this paper was undertaken. This work was supported by EPSRC grant no EP/K032208/1. While working on this project the second author was supported by DFG Schwerpunktprogramm 1786 “Homotopy Theory and Algebraic Geometry” (project ID SCHW 860/1-1). He is currently an associate member of the SFB 1085 Higher Invariants. 


\section{Cyclification of global spaces}\label{sec:cyc}

Our first goal is to construct the cyclification functor. We begin with some recollections on global spaces.

\begin{mydef}
Let $\Glo$ be the \emph{global orbit category} whose objects are compact Lie groups $G$ and whose morphism spaces are given by $\Map_{\Glo}(H,G) = |\Hom(H,G) // G|$, the geometric realisation of the action groupoid of $G$ acting by conjugation on the space of homomorphisms of compact Lie groups. Let $\Spc_\gl = \Fun(\Glo^\op, \Spc)$ be the category of \emph{global spaces}. We write the global space represented by a group $G$ by $\B G$. Let $\Glo_{\ab\fin}$ be the subcategory of $\Glo$ spanned by those compact Lie groups $G$ which are finite and abelian; this is equivalent to the subcategory of spaces spanned by those spaces of the form $BA$ for a finite abelian group $A$. Let $\Spc_{\ab\fin} = \Fun(\Glo^\op_{\ab\fin},\Spc)$ denote the subcategory of $\Spc_\gl$ generated under colimits by representables from $\Glo_{\ab\fin}$.
\end{mydef}

The category $\Spc_\gl$ was first introduced (as a topological category) by Gepner--Henriques in \cite{gepnerheriques}. Schwede gives an alternative construction in \cite[\textsection1]{s} and also compares these definitions in \cite{schwedeglocomparison}. Gepner--Meier give an $\infty$-categorical treatment in \cite[\textsection2]{davidandlennart}. On the other hand, Lurie uses $\Glo_{\ab\fin}$ and $\Spc_{\ab\fin}$ throughout \cite{ec3}, where he denotes them by $\mathscr{T}$ and $\mathcal{OS}$, respectively.

\begin{mydef}
Consider the functor $\mathrm{ev}_{\ast}\colon \Spc_{\gl}\rightarrow \Spc$ sending a global space $X$ to $X(\ast)$, the evaluation of $X$ on $\ast \in \Glo$. This functor admits a right adjoint denoted by $\R\colon \Spc\rightarrow \Spc_{\gl}$. Explicitly, $\R(X)$ is given by the assignment $\B H\mapsto \Map(B H,X)$. For a $1$-truncated compact Lie group $G$, we have an equivalence $\R(BG) \simeq \B G$ by \cite[Th.1.2]{rezktruncatedliegroups}, where $BG = \ast/G$ in $\Spc$. In particular, this holds for all abelian compact Lie groups $G$.
\end{mydef}

Now we can define the cyclification functor.

\begin{mydef}\label{cyclificationdefinition}
We write $\mathbf{T}$ for the compact Lie group $U(1)$. Let $\Cyc \colon \Spc_{\ab\fin}\to \Spc_{\gl}$ be the left Kan extension along the Yoneda embedding $\Glo_{\ab\fin}\to \Spc_{\gl}$ of the composite
\[\Glo_{\ab\fin} \to \Spc \xrightarrow{L(-)/\TT} \Spc \xrightarrow{\R} \Spc_{\gl}\]
where the first functor is the canonical inclusion and the second is the quotient of the free loop space functor $L(-)$ by the natural $\TT$-action coming from the rotation of loops.
\end{mydef}

\begin{example}\label{valueofcyconapoint}
The value of $\Cyc(\ast)$ is equivalent to $\B\TT$:
\[\R(L(\ast)/\TT) \simeq \R(\ast /\TT) \simeq \R(B\TT) \simeq \B \TT \in \Spc_\gl\]
\end{example}

\begin{example}\label{valueofcyconacyclic}
More generally, one can calculate $\Cyc$ on $\B A$, where $A$ is a finite abelian group. First, we recall the identification $L( B A) \simeq \coprod_{a\in A} BA$. Under this equivalence, the $\TT$-action on $L(B A)$ factors through a $\TT$-action on each coproduct factor. By inspection, it agrees on the factor associated with the element $a\in A$ with the restriction of the canonical action of $BA$ on itself (using that $BA$ is a group) along the map $B\Z\to BA$ induced by the element $a\in A$. In this way we obtain an equivalence $L(BA)/\TT \simeq \coprod_{a\in A} X_a$, where $X_a$ sits in the following pullback squares
\[\begin{tikzcd}
	BA & {X_a} & \ast \\
	\ast & {B^2\Z} & {B^2A.}
	\arrow[from=1-1, to=1-2]
	\arrow[from=1-1, to=2-1]
	\arrow["\lrcorner"{anchor=center, pos=0.125}, draw=none, from=1-1, to=2-2]
	\arrow[from=1-2, to=1-3]
	\arrow[from=1-2, to=2-2]
	\arrow["\lrcorner"{anchor=center, pos=0.125}, draw=none, from=1-2, to=2-3]
	\arrow[from=1-3, to=2-3]
	\arrow[from=2-1, to=2-2]
	\arrow["a", from=2-2, to=2-3]
\end{tikzcd}\]
Writing $n$ for the order of $a\in A$, the map $a\colon B^2\Z \to B^2 A$ factors through $B^2 \Z/n$. Now we can extend the diagram above further to
\[\begin{tikzcd}
	BA & {X_a} & {BA/\langle a\rangle} & \ast \\
	\ast & {B^2\Z} & {B^2\Z/n} & {B^2A} \\
	&& \ast & {B^2 A/\langle a\rangle.}
	\arrow[from=1-1, to=1-2]
	\arrow[from=1-1, to=2-1]
	\arrow["\lrcorner"{anchor=center, pos=0.125}, draw=none, from=1-1, to=2-2]
	\arrow[from=1-2, to=1-3]
	\arrow[from=1-2, to=2-2]
	\arrow["\lrcorner"{anchor=center, pos=0.125}, draw=none, from=1-2, to=2-3]
	\arrow[from=1-3, to=1-4]
	\arrow[from=1-3, to=2-3]
	\arrow["\lrcorner"{anchor=center, pos=0.125}, draw=none, from=1-3, to=2-4]
	\arrow[from=1-4, to=2-4]
	\arrow[from=2-1, to=2-2]
	\arrow[from=2-2, to=2-3]
	\arrow[from=2-3, to=2-4]
	\arrow[from=2-3, to=3-3]
	\arrow["\lrcorner"{anchor=center, pos=0.125}, draw=none, from=2-3, to=3-4]
	\arrow[from=2-4, to=3-4]
	\arrow[from=3-3, to=3-4]
\end{tikzcd}\]
As $\Z/n$ acts freely on $A$, the map $BA/\langle a\rangle \to B^2 \Z/n$ factors through the point, and so we obtain the diagram
\[\begin{tikzcd}
	{X_a} & {B^2 \Z} & {B^2\Z} \\
	{B A/\langle a\rangle} & \ast & {B^2 \Z/n,}
	\arrow[from=1-1, to=1-2]
	\arrow[from=1-1, to=2-1]
	\arrow["\lrcorner"{anchor=center, pos=0.125}, draw=none, from=1-1, to=2-2]
	\arrow[from=1-2, to=1-3]
	\arrow[from=1-2, to=2-2]
	\arrow["\lrcorner"{anchor=center, pos=0.125}, draw=none, from=1-2, to=2-3]
	\arrow[from=1-3, to=2-3]
	\arrow[from=2-1, to=2-2]
	\arrow[from=2-2, to=2-3]
\end{tikzcd}\]
giving an equivalence $X_a \simeq B^2 \Z \times B A/\langle a\rangle \simeq B(\TT \times A/\langle a\rangle)$. As $\R$ sends classifying spaces of $1$-truncated compact Lie groups to representables, we find that 
\[
\Cyc(\B A) \simeq \coprod_{a\in A} \B (\TT\times A/\langle a\rangle).
\]
In particular, the functor $\Cyc\colon \Spc_{\ab\fin} \to \Spc_\gl$ restricted to $\Glo_{\ab\fin}$ factors through the subcategory $\Glo^{\sqcup}_\ab \subseteq \Spc_{\gl}$ of those global spaces given by finite coproducts of representables $\B G$ for abelian $G$:
\begin{equation}\label{factorisationofcyc}
   \Cyc\colon \Glo_{\ab\fin} \to \Glo_\ab^\sqcup
\end{equation}
\end{example}

\begin{remark}[GRH-extended inertia orbifold]\label{cycislambda}
We may compare this to the \emph{extended inertia orbifold} of Ganter, Rezk, and Huan. Again, let $A$ be a finite abelian group. Recall that $\Lambda_{\mathrm{ext}}^{\mathrm{GRH}}(\B A)$ is defined to be $\coprod_{a\in A} \B \Lambda_{a/A}$, where $\Lambda_{a/A}$ is a compact Lie group defined by the quotient
\[
0\rightarrow \Z \rightarrow A\times \R \rightarrow \Lambda_{a/A}\rightarrow 0.
\] Alternatively, if we write $n$ for the order of $A$, then we also have a short exact sequence
\[0 \to \Z/n \to A\times \TT \to \Lambda_{a/A} \to 0.\]
Suppose $f\colon A\to A'$ is a group homomorphism such that $f(a)= a'$. Write $m = |a|/|a'|$ for the quotient of the order of $a$ by the order of $a'$. This induces a group homomorphism $\Lambda_{a/A}$ by the quotient of the map 
\[
(f,m\cdot (-))\colon A\times \TT \to A'\times \TT.
\]
Applying this to each component gives the functoriality of $\Lambda^{\mathrm{GRH}}_{\mathrm{ext}}$ in group homomorphisms. To connect to the previous discussion, we note that $\Lambda_{a/A}$ sits in the canonically split short exact sequence
\[0 \to \TT \xrightarrow{t\mapsto [e,t]} \Lambda_{a/A} \xrightarrow{[x,t]\mapsto [x]} A/\langle a \rangle \to 0,\]
with retraction given by $\Lambda_{a/A} \to \TT$ sending $[x,t]$ to $t$. In other words, we have a natural identification
\begin{equation}\label{splittingoflambda}    
\Lambda_{a/A} \simeq \TT \times A/\langle a\rangle,  
\end{equation}
and so an identification 
\[
\Lambda_{\mathrm{ext}}^{\mathrm{GHR}}(\B A) \simeq \Cyc(\B A),
\]
which is easily shown to be natural in maps of groups. It is less clear how to obtain an identification that is natural in maps of groupoids. In fact, we are not sure how to apriori make the left-hand side of the equivalence a functor in groupoids. This is why we prefer our more invariant definition of cyclification.
\end{remark}

\begin{remark}[Cyclification of Sati--Schreiber]
Our choice of notation and the definition of $\Cyc$ is inspired by \cite{satischreiber}, who define the cyclification of arbitrary sheaves on smooth manifolds and compare this to the extended inertia orbifold definition of Ganter, Rezk and Huan for ``good orbifolds''. We have nevertheless chosen to give a different perspective on cyclification, which we find simpler. Namely, our definition exploits the fact that $\Glo_{\ab}$ and the category of classifying spaces of compact abelian Lie groups are equivalent, allowing us to directly define cyclification in spaces. A definition for more general global spaces, such as those of the form $\B G$ for a general compact Lie group $G$, is possible by following the definitions of \cite{satischreiber} and using \cite{clough2024globalspaceshomotopytheory}.
\end{remark}

The cyclification functor behaves well on products of finite abelian groups.

\begin{prop}\label{productsovercircle}
Let $A,A'$ be finite abelian groups. Then the natural map
\[\Cyc(\B A\times \B A') \to \Cyc(\B A)\times_{\BT} \Cyc(\B A')\]
is an equivalence of global spaces.
\end{prop}

\begin{proof}
First, we have the natural identification $L(B A\times B A') \simeq L(BA) \times L(BA')$. Both sides have obvious $\TT \times \TT$-actions, and we may compute that 
\[L(BA \times BA')/\TT \times \TT \simeq L(B A)/\TT \times L(BA')/\TT.\]
By descent, the quotient $L(B A\times B A')/
\TT$ therefore sits in a pullback square
\[\begin{tikzcd}
	{L(BA\times BA')/\TT} & {L(BA)/\TT \times L(BA')/\TT} \\
	B\TT & {B\TT\times B\TT.}
	\arrow[from=1-1, to=1-2]
	\arrow[from=1-1, to=2-1]
	\arrow[from=1-2, to=2-2]
	\arrow[""{name=0, anchor=center, inner sep=0}, "\Delta", from=2-1, to=2-2]
	\arrow["\lrcorner"{anchor=center, pos=0.125}, draw=none, from=1-1, to=0]
\end{tikzcd}\]
This implies that $L(BA\times BA')/\TT \simeq L(BA)/\TT \times_{B\TT} L(BA')/\TT$. Finally, because $\R$, as a right adjoint, preserves pullbacks, we are done.
\end{proof}


\section{Cyclification of oriented equivariant cohomology theories}\label{sec:cyclifyell}

The main result of this section is \Cref{emergenceofpdivgroup}, which states that precomposing an equivariant cohomology theory with the cyclification functor of the previous section produces an equivariant cohomology theory. More specifically, we will be using \cite{davidandlennart} and \cite{ec3} to produce our equivariant cohomology theories, so we begin by recalling some of these notions.\\

Fix an affine stack $\s = \Spec A$ for this section. Recall that an abelian group object in a category $\C$ with finite products is a product preserving functor $F\colon \Lat^\op \to \C$ where $\Lat$ is the category of lattices, so finitely generated free abelian groups.\\

Let $\G$ a preoriented abelian stack over $\s$, so an object $\G$ of $\Ab(\Stk_{/\s})$ equipped with a morphism of abelian group objects $e\colon B\TT\to \G(\s)$. By applying \cite[Con.3.13]{davidandlennart} and then extending via finite coproducts, we may extend $\G$ to a functor
\begin{equation}\label{ellcohom}\Ell_{\G/\s}\colon \Glo_{\ab}^\sqcup \to \Shv(\s), \qquad \B H \mapsto \HOM_\s(\underline{H}^\vee, \G),\end{equation}
where $H^\vee = \Hom(H, \TT)$ denotes the Pontryagin dual. In particular, its value on $\B\TT$ recovers $\G$ itself
\begin{equation}\label{ellonbt}
    \Ell_{\G/\s}(\B\TT) \simeq \HOM_\s(\underline{\TT}^\vee, \G) \simeq \HOM_\s(\underline{\Z}, \G)\simeq \G ,  
\end{equation}
and its value on cyclic groups recovers the torsion of $\G$
\begin{equation}\label{ellonbcn}
    \Ell_{\G/\s}(\B \Z/n) \simeq \HOM_\s(\underline{\Z/n}^\vee, \G) \simeq  \G[n] .  
\end{equation}

\begin{remark}\label{rmk:elliptic_cohom_always_a_little_exact}
As observed in the proof of \cite[Pr.3.15]{davidandlennart}, see also \cite[Pr.2.16]{davidandlennart2}, the functor $\Ell_{\G/\s}(-)\colon \Glo_{\ab}\to \Stk_{/\s}$ preserves finite products and pullbacks in which at least one of the maps is induced by a surjective map of groups. Because coproducts are disjoint in $\Stk_{/\s}$, the extension to $\Glo_{\ab}^{\sqcup}$ again preserves pullbacks for which the same condition holds on each component.
\end{remark}

We work with a subcategory of preoriented abelian stacks for which the torsion of $\G$ is well-behaved. 

\begin{mydef}\label{defofgeometricstacks}
    Let $\G$ be an abelian stack over $\s$. We say that $\G$ is \emph{geometric} if the formal completion of $\G$ is a formal group, the map $[n]\colon \G\to \G$ is finite flat for all integers $n\geq 1$, and $\G$ is almost of finite presentation over $\x$. From the proof of \cite[Pr.6.7.1]{ec1}, it follows that $\G[\P^\infty]$ is a $\P$-divisible group over $\s$, in the sense of \Cref{df:preorientedpdivgp}.

\end{mydef}

Examples of such geometric abelian stacks are abelian varieties (\cite[Pr.6.7.1]{ec1} and \cite[Pr.7.1.2]{ec2}) and tori (\cite{geometricnorms}), the latter being an abelian stack which fpqc locally looks like $\G_m$. In \cite[Df.3.6]{elltempcomp}, the first-named author defines and discusses geometric abelian sheaves, a slightly more general, but related notion.\\

The next statement follows almost straight from this.

\begin{prop}\label{ellcohomfactoring}
    Let $\G$ be a geometric preoriented abelian $\s$-stack. Write $\Glo_\ab^\sqcup$ for the subcategory of $\Spc_\gl$ generated by finite coproducts of representables $\B H$ where $H$ is an abelian compact Lie group. Then the restriction of $\Ell_{\G/\s}$ of (\ref{ellcohom}) factors through $\Stk_{/\s}^\afp$, the category of stacks over $\s$ which are almost of finite presentation:
    \[\Ell_{\G/\s}\colon \Glo^\sqcup_\ab \to \Stk^\afp_{/\s}.\]
\end{prop}

\begin{proof}
This essentially follows from (\ref{ellonbt}), (\ref{ellonbcn}), and the classification of abelian compact Lie groups; \cite[Rmk.1.10]{elltempcomp}.
\end{proof}

Now that we have identified our desired equivariant cohomology theories $\Ell_{\G/\s}$, we can start to cyclify them.

\begin{mydef}\label{df:cyclicdcationofellcohom}
    Let $\s$ be a stack and $\G$ be a geometric preoriented abelian $\s$-stack. The factorisations (\ref{factorisationofcyc}) and (\ref{ellcohomfactoring}) yield the composition
    \[\Glo_{\ab\fin} \xrightarrow{\Cyc} \Glo_{\ab}^\sqcup \xrightarrow{\Ell_{\G/\s}} \Stk_{/\s}\]
    which sends the terminal object $\ast$ of $\Glo_{\ab\fin}$ to $\Ell_{\G/\s}(\Cyc(\ast))\simeq \Ell_{\G/\s}(\B\TT) \simeq \G$ by \Cref{valueofcyconapoint} and (\ref{ellonbt}), respectively. Let us write $\CycEll_{\G/\s}$ for the unique factorisation of the above functor through the forgetful functor $\Stk_{/\G}\to \Stk_{/\s}$
    \[\CycEll_{\G/\s} \colon \Glo_{\ab\fin} \to \Stk_{/\G}\]
    and call it the \emph{cyclification} of $\Ell_{\G/\s}$.
\end{mydef}

Our next goal is to show that if $\G$ is geometric, then the cyclification of $\Ell_{\G/\s}$ is a preoriented $\P$-divisible group over $\G$. First, we recall this notion.

\begin{mydef}\label{df:preorientedpdivgp}
    Let $\s$ be a stack. A \emph{preoriented $\P$-divisible group} over $\s$ is a functor $F\colon \Glo_{\ab\fin} \to \Stk_{/\s}$ satisfying the following conditions:
    \begin{enumerate}
        \item The functor $F$ preserves finite products. In particular, there is a natural equivalence $F(\ast) \simeq \s$.
        \item For every short exact sequence of finite abelian groups $0\to A'\to A\to A''\to 0$, the diagram of $\s$-stacks
        \[\begin{tikzcd}
            {F(A')}\ar[r]\ar[d]    &   {F(0)}\ar[d]   \\
            {F(A)}\ar[r]           &   {F(A'')}
        \end{tikzcd}\]
        is Cartesian. Moreover, the horizontal maps are finite flat of positive degree.
    \end{enumerate}
\end{mydef}

\begin{remark}[Over an affine base]\label{reductiontofinitecase}
    In \cite[Th.3.5.5]{ec3}, Lurie's shows that if $\s = \Spec R$ is affine, then a preoriented $\P$-divisible group $\G$ as defined above is precisely the data of a preoriented $\P$-divisible group $\G$ over $R$ \`{a} la \cite[\textsection 6.5]{ec1} together with a preorientation $e$ of $\G$ \`{a} la \cite[\textsection2.6]{ec3}.
\end{remark}

This leads to the main theorem of this section.

\begin{theorem}\label{emergenceofpdivgroup}
    Let $\s$ be a stack and $\G$ be an geometric abelian $\s$-stack. Then \[
    \CycEll_{\G/\s}\colon \Glo_{\ab\fin}\to \Stk_{/\G}
    \] defines a preoriented $\P$-divisible group over $\G$.
\end{theorem}




\begin{proof}[Proof of \Cref{emergenceofpdivgroup}]
    For simplicity, let us write $\EE = \Ell_{\G/\s}$ and $\C = \CycEll_{\G/\s}$. We now check that $\C$ satisfies the conditions of \Cref{df:preorientedpdivgp}. For part 1, we have $\C(\ast)\simeq \EE(\B\TT)\simeq \G$ by construction; see \Cref{df:cyclicdcationofellcohom}. For two finite abelian groups $A,A'$, we consider the composition
    \[\EE(\Cyc(\B A \times \B A')) \xrightarrow{f} \EE\left(    \Cyc(\B A) \underset{\Cyc(\ast)}{\times} \Cyc(\B A')    \right) \xrightarrow{g} \EE(\Cyc(\B A)) \underset{\EE(\Cyc(\ast))}{\times} \EE(\Cyc(\B A')).\]
    The map $f$ is an equivalence from \Cref{productsovercircle}. On the other hand, the map $\Cyc(\B A)\to \Cyc(\ast)$ is surjective on each component for all $A$, and so the map $g$ is an equivalence by \Cref{rmk:elliptic_cohom_always_a_little_exact}. This shows part 1.
    
    For part 2, consider a short exact sequence of abelian groups
    \[0 \to A' \to A \xrightarrow{p} A'' \to 0.\]
    Let us write $A_a=A/\langle a\rangle$ for the quotient of $A$ by the subgroup generated by $a$. Using that $\C$ preserves products, the natural identification of abelian groups (\ref{splittingoflambda}), and the fact that $\EE$ preserves products, we naturally identify the two commutative squares of $\s$-stacks
    \[\begin{tikzcd}
        {\C(\B A')}\ar[r]\ar[d] &   {\C(\ast)}\ar[d]   \\
        {\C(\B A)}\ar[r]        &   {\C(\B A'')}
    \end{tikzcd}\qquad\qquad 
    \begin{tikzcd}
        {\coprod_{a'\in A'} \EE(\B A'_{a'})\underset{\s}{\times} \G}\ar[r]\ar[d]    &   {\G}\ar[d] \\
        {\coprod_{a\in A} \EE(\B A_a)\underset{\s}{\times} \G}\ar[r]              &  {\coprod_{a''\in A''} \EE(\B A''_{a''})\underset{\s}{\times} \G.}
    \end{tikzcd}\]
    To show that the diagrams above are Cartesian, we work over each component $a'' \in A''$ in the lower-right corner of the right square. The right-hand vertical map hits the component associated with $a'' = 0$, while the left-hand vertical map also only hits components such that $p(a) = 0$. Therefore over any $a''\neq 0$, the map is trivially a pullback, since the empty stack pulls back to the empty stack along any map of stacks. Therefore, it suffices to work over $a'' = 0$, in other words, it suffices to show that the diagram 
    \[
    \begin{tikzcd}
        {\coprod_{a'\in A'} \EE(\B A'_{a'})\underset{\s}{\times} \G}\ar[r]\ar[d]    &   {\G}\ar[d] \\
        {\coprod_{a'\in A' \subseteq A} \EE(\B A_{a'})\underset{\s}{\times} \G}\ar[r]              &  \EE(\B A'')\underset{\s}{\times} \G
    \end{tikzcd}\]
    is Cartesian. We may factor this through the fold map
    \[
    \begin{tikzcd}
	{\coprod_{a'\in A'} \EE(\B A'_{a'})\underset{\s}{\times}\G} & {\coprod_{a'\in A'}\G} & \G \\
	{\coprod_{a'\in A'} \EE(\B A_{a'})\underset{\s}{\times}\G} & {\coprod_{a'\in A'}\EE(\B A'')\underset{\s}{\times}\G} & {\EE(\B A'')\underset{\s}{\times}\G.}
	\arrow[from=1-1, to=1-2]
	\arrow[from=1-1, to=2-1]
	\arrow["\nabla", from=1-2, to=1-3]
	\arrow[from=1-2, to=2-2]
	\arrow[from=1-3, to=2-3]
	\arrow[from=2-1, to=2-2]
	\arrow["\nabla", from=2-2, to=2-3]
    \end{tikzcd}\]
    Since the right-hand square is a pullback, it suffices to show that each component of the left-hand square is a pullback. In other words that 
    \[\begin{tikzcd}
	{\EE(\B A'_{a'})\underset{\s}{\times}\G} & \G \\
	{\EE(\B A_{a'})\underset{\s}{\times}\G} & {\EE(\B A'')\underset{\s}{\times}\G}
	\arrow[from=1-1, to=1-2]
	\arrow[from=1-1, to=2-1]
	\arrow[from=1-2, to=2-2]
	\arrow[from=2-1, to=2-2]
\end{tikzcd}\]
is a pullback for all $a'\in A'$. Since all of the maps in the diagram above are induced by the functoriality of pullbacks, we may compute the pullback in each variable separately. These are 
\[\begin{tikzcd}
	{\EE(\B A'_{a'})} & \s \\
	{\EE(\B A_{a'})} & {\EE(\B A'')}
	\arrow[from=1-1, to=1-2]
	\arrow[from=1-1, to=2-1]
	\arrow[from=1-2, to=2-2]
	\arrow[from=2-1, to=2-2]
\end{tikzcd} \qquad\text{and}\qquad \begin{tikzcd}
	\G & \G \\
	\G & \G
	\arrow["{-\cdot|a'|}", from=1-1, to=1-2]
	\arrow["\id"', from=1-1, to=2-1]
	\arrow["\id", from=1-2, to=2-2]
	\arrow["{-\cdot |a'|}", from=2-1, to=2-2]
\end{tikzcd}\]
respectively. The second square is clearly a pullback. The first is a pullback, as it is given by applying the $\P$-divisible group $\EE$ to the short exact sequence
\[
0 \to A'/\langle a'\rangle \to A/\langle a'\rangle \to A'' \to 0
\]
of abelian groups.
\end{proof}

Although the assignment sending a geometric abelian $\x$-stack $\G$ to the preoriented $\P$-divisible group $\CycEll_{\G/\x}$ is a way of producing a new interesting equivariant cohomology from $\G$, it does not increase the chromatic height of $\G$, only the $\P$-divisible height. For example, as we will see below, cyclifying $\G_m/\KU$ produces something akin to the torsion of an \emph{ordinary} elliptic curve, meaning an elliptic curve with formal height $1$.


\section{A spectral lift of the Katz--Mazur \texorpdfstring{$\P$}{P}-divisible group}\label{kmsection}

The Katz--Mazur group scheme $\T_\km^\heartsuit$ is a $\P$-divisible group defined over $\G_{m,\Z}$ such that its associated torsion points base-change over $\Z\llpar q \rrpar$ to the $\P$-divisible group associated with the Tate curve; see \cite[\textsection8.7]{km}. As we will show now, the previous section provides a functorial lift of $\T_\km^\heartsuit$ to a group object over $\G_{m,\KU}$.\\

Let $\G_m$ be the multiplicative group scheme over the sphere from \cite[\textsection1.6.3]{ec2} whose underlying spectral scheme is given by the spectrum of $\Sph[q^{\pm}]=\Sph[\N][q^{-1}]=\Sph[\Z]$, the localisation of the suspension spectrum of the natural numbers at the generator $q$. By Snaith's theorem, see \cite[\textsection6.5]{ec2}, $\KU$ is the orientation classifier of $\G_m$, and so $\G_m$ has a canonical orientation after base-change to $\KU$. The machinery of \cite{davidandlennart} then gives us the functor (\ref{ellcohom}), which now takes the form
\[\Ell_{\G_m/\KU} \colon \Glo_{\ab}^{\sqcup} \to \Shv(\Spec \KU);\]
this is described in detail in \cite[\textsection4]{davidandlennart}, where it is also shown to agree with the cohomology theory on global spaces given by equivariant $K$-theory $\gKU$. Per construction, the value of $\Ell_{\G_m/\KU}$ on $\B\TT$ is the affine group scheme $\G_{m}$ over $\KU$. From its description over $\Sph$ above, this spectral scheme is affine and its global sections are given by the $\E_\infty$-ring $\KU[q^\pm]$.\\

We can now cyclify $\Ell_{\G_m/\KU}$ as done in the previous section and produce an {oriented} $\P$-divisible group.

\begin{theorem}\label{spectralliftofkm}
The preoriented $\P$-divisible group $\CycEll_{\G_m/\KU}$ over $\Spec \KU[q^\pm]$ is an oriented $\P$-divisible group, see \cite[Df.2.6.12]{ec3}, now denoted as $\T_\km$, with underlying classical $\P$-divisible group the Katz--Mazur group scheme $\T_\km^\heartsuit$ over $\Z[q^\pm]$.
\end{theorem}

\begin{remark}[Orbifold quasi-elliptic cohomology]\label{rmk:compwithgrh}
Taking global sections and passing to homotopy groups we obtain a global cohomology theory $\pi_{-*}\Gamma \CycEll_{\G_m/\KU}$ valued in graded commutative rings. This definition has appeared in the literature before. See for example \cite[Df.4.5]{huanquasielliptic}, where it is called \emph{orbifold quasi-elliptic cohomology}. The agreement of our two definitions follows immediately from the agreement of $\Cyc(-)$ and $\Lambda^{\mathrm{GRH}}_{\mathrm{ext}}(-)$ discussed in \Cref{cycislambda}. We explore and refine this further in \Cref{ssec:qell}.
\end{remark}

\begin{proof}[Proof of \Cref{spectralliftofkm}]
The identification of the underlying classical $\P$-divisible group of $\T_\km$ goes back to \cite{rezkquasielliptic} at least. One proceeds by using the natural identification of $\gKU^0(\B(-))$ with the functor $\mathrm{Rep}(-)$ of complex representation rings. This allows one to explicitly compute all of the rings $\gKU^0(\Cyc(\B G))$ as well as the functoriality in group homomorphisms on finite abelian groups $G$ using \Cref{valueofcyconacyclic,productsovercircle}. The resulting formulas agree precisely with those appearing in the definition of $\T_\km^\heartsuit$ by inspection. See for example \cite[Rmk.3.13]{huanquasielliptic}, where Huan shows that the spectrum of $\pi_{0} \gKU(\Cyc(\B \Z/N)) \simeq \gKU^0(\Lambda^{\mathrm{GHN}}_{\mathrm{ext}} \B \Z/N)$ is naturally equivalent to the $N$-torsion of $\T_\km^\heartsuit$ as group schemes.\\
       
Now onto the orientability of $\CycEll_{\G_m/\KU} = \T_\km$. Let us also write $\A = \KU[q^\pm]$ for brevity. We already know from \Cref{emergenceofpdivgroup} that $\T_\km$ is preoriented, so we are left to show that the natural preorientation on $\T_\km$ is actually an orientation. As the definition of a preorientation being an orientation only depends on the $p$-completion of the base stack for each prime $p$, we will now work over the $p$-completion of $\A$, denoted by $\widehat{\A}$, and write $\widehat{\T}_\km$ for the associated $p$-divisible group over $\widehat{\A}$. Notice that $\widehat{\A}$ is complex-periodic à la \cite[Df.4.1.8]{ec2}. By \cite[Pr.2.2.1]{ec3}, a preorientation on $\widehat{\T}_\km$ is equivalent to a preorientation on its identity component $\widehat{\T}_\km^\circ$, and by \cite[Pr.4.3.21]{ec2}, preorientations on $\widehat{\T}_\km^\circ$ are precisely maps of formal groups $\widehat{\G}^\QQ_{\widehat{\A}} \to \widehat{\T}_\km^\circ$ over $\widehat{\A}$, where $\widehat{\G}^\QQ_{\widehat{\A}}$ is the \emph{Quillen formal group} of the complex-periodic $\E_\infty$-ring $\widehat{\A}$; see \cite[\textsection4.1]{ec2}. To show that our preorientation $e\colon \widehat{\G}^\QQ_{\widehat{\A}} \to \widehat{\T}_\km^\circ$ is an equivalence, it suffices to show this over $\pi_0 {\widehat{\A}} \simeq \Z[q^\pm]^\wedge_p$, as formal groups over $\E_\infty$-rings are flat over the base. This is now clear. Indeed, the complex orientation on $\widehat{\A}$ comes from $\KU$, hence $\widehat{\G}^\QQ_{\widehat{\A}}$ can be identified with the formal multiplicative group $\widehat{\G}_m$ over $\widehat{\A}$. Moreover, the identity component of $\widehat{\T}_\km$ over $\pi_0 \widehat{\A}$ can also be identified with $\widehat{\G}_m$ as $\T_\km^\heartsuit$ is the classical Katz--Mazur group scheme over $\Z[q^\pm]$ which sits in an exact sequence
    \[0 \to \G_m \to \T_\km^\heartsuit \to \underline{\Q/\Z} \to 0\]
    of fppf sheaves over $\Z[q^\pm]$ by \cite[(8.7.2.1)]{km}.
\end{proof}

\begin{mydef}\label{kmschemedef}
    Write $\T_\km$ for the oriented $\P$-divisible group of \Cref{spectralliftofkm}. In other words, $\T_\km$ is the functor 
    \[
    \T_\km = \CycEll_{\G_m/\KU} \colon \Glo_{\ab\fin} \to \Stk_{/\Spec \KU[q^{\pm}]}
    \] of \Cref{emergenceofpdivgroup}, which happens to be an oriented $\P$-divisible group. By \Cref{reductiontofinitecase}, this is a $\P$-divisible group $\T_\km$ over the $\E_\infty$-ring $\KU[q^{\pm}]$ together with a choice of orientation $e: \Sigma(\Q/\Z) \to \T_\km(\KU\llpar q \rrpar)$.
\end{mydef}

In \Cref{ssec:qell}, we will further study the tempered cohomology theory associated with $\T_\km$. Before we do that, we would like to show that $\T_\km$ is actual the torsion of an elliptic curve when base-changed to $\KU\llpar q \rrpar$.


\section{The derived Tate curve}\label{tatecurvesection}

The $\P$-divisible group $\T_\km$ over $\G_{m,\KU}$ is not yet the $\P$-divisible group of an elliptic curve. Indeed, even classically we need to base-change from $\Z[q^\pm]$ to $\Z\llpar q\rrpar$ in order to introduce the classical Tate curve $\T^\heartsuit$. We would like to do the same spectrally. Inverting $q$ on both sides of the completion $\KU[q]\to \KU\llbracket q\rrbracket$ induces a map
\[    \KU[q^\pm]\to \KU\llpar q\rrpar,\]
and base-changing $\T_\km$ along this map produces an oriented $\P$-divisible group over $\KU\llpar q \rrpar$. The following is a slight refinement of \Cref{existenceoftatecurve}.

\begin{theorem}\label{liftoftate}
     There exists an oriented elliptic curve $\T$ over $\KU\llpar q\rrpar$ such that:
     \begin{enumerate}
    \item its underlying classical elliptic curve is uniquely isomorphic to the classical Tate curve $\T^\heartsuit$ over $\Z\llpar q\rrpar$,
    \item its associated formal group is the formal multiplicative group scheme $\widehat{\G}_{m}$ over $\KU\llpar q \rrpar$,
    \item the associated map $\TMF \to \KU\llpar q \rrpar$ induces the classical meromorphic $q$-expansion map of (\ref{qexpmap}) over $\Q$ on rational homotopy groups, and 
    \item its associated oriented $\P$-divisible group is the base-change of $\T_\km$ from $\KU[q^\pm]$ to $\KU\llpar q\rrpar$.
\end{enumerate}
     Moreover, such a $\T$ is uniquely determined up to homotopy as an oriented elliptic curve over $\KU\llpar q \rrpar$ by:
     \begin{itemize}
         \item Part 3 above,
         \item An isomorphism between the associated $p$-divisible group over $\KU\llpar q \rrpar^\wedge_p$ and that of $\T_\km$ base-changed to $\KU\llpar q \rrpar^\wedge_p$ for each prime $p$,
         \item and an isomorphism between its underlying classical $p$-divisible group base-changed and the $p$-divisible group of the classical Tate curve over $\Z\llpar q\rrpar^\wedge_p$.
     \end{itemize} 
\end{theorem}

This uniqueness statement, a by-product of our construction of $\T$, would also be useful to identify various other potential definitions of a derived Tate curve over $\KU\llpar q \rrpar$. Another potential uniqueness statement could be made in terms of the $\KU$-homology of the associated map $\TMF \to \KU\llpar q \rrpar$, proven using Goerss--Hopkins obstruction theory similar to \cite[Th.A]{realspectra}. One could hopefully use the uniqueness of $\T$ in \Cref{liftoftate} to show that our construction of $\T$ agrees with Smit's forthcoming work of the Tate curve over $\Sph\llbracket q \rrbracket$ base-changed to $\KU\llpar q \rrpar$.\\

Our proof of \Cref{liftoftate} requires two ingredients: an oriented variant of Lurie's spectral Serre--Tate theorem (\Cref{orientedspectralserretate}) and a corepresentability property of $\TMF$ (\Cref{corepresentabilityoftmfinnicecases}).\\

Recall that we write $\Ell^\ori(R)$ for the coCartesian unstraightening of the functor $\Ell(R) \to \Spc$ sending an elliptic curve over an $\E_\infty$-ring $R$ to its space of orientations; see \cite[Not.7.2.4]{ec2}. Similarly, we write $\BTtwoor(R)$ for coCartesian unstraightening of the functor $\BTtwo(R) \to \Spc$ sending a $p$-divisible group of height $2$ over an $\E_\infty$-ring $R$ to its space of orientations.

\begin{theorem}[Oriented spectral Serre--Tate theorem]\label{orientedspectralserretate}
    Let $p$ be a prime and $R$ be a $p$-complete $\E_\infty$-ring, then the diagram of categories
    \[\begin{tikzcd}
        {\Ell^\ori(R)}\ar[r]\ar[d, "{[p^\infty]}"]   &   {\Ell(\pi_0 R)}\ar[d, "{[p^\infty]}"]   \\
        {\BTtwoor(R)}\ar[r]   &   {\BTtwo(\pi_0 R)}
    \end{tikzcd}\]
    is Cartesian.
\end{theorem}

This will follow from the following slight reformulation of \cite[Th.7.0.1]{ec1} together with an identification of spaces of orientations of elliptic curves and their associated $p$-divisible groups over $p$-complete $\E_\infty$-rings. Recall that for an $\E_\infty$-ring $R$, we write $\AVar_g(R)$ for the category of \emph{(strict) abelian varieties of dimension $g$} (\cite[Df.1.5.1]{ec1}) and $\BTn(R)$ for the category of \emph{$p$-divisible groups of height $n$} (\cite[Df.6.5.1]{ec1}). 

\begin{theorem}[Spectral Serre--Tate theorem]\label{spectralserretate}
    Let $p$ be a prime and $R$ be a $p$-complete $\E_\infty$-ring, then the diagram of categories
    \[\begin{tikzcd}
        {\AVar_g(R)}\ar[r]\ar[d, "{[p^\infty]}"]   &   {\AVar_g(\pi_0 R)}\ar[d, "{[p^\infty]}"]   \\
        {\BTg(R)}\ar[r]   &   {\BTg(\pi_0 R)}
    \end{tikzcd}\]
    is Cartesian.
\end{theorem}

The only difference between the above theorem and \cite[Th.7.0.1 \& Pr.7.4.1]{ec1} is that we do not require there to exist a morphism of $\E_\infty$-rings $R\to \pi_0 R$, which never exists if $R$ is nonzero and complex-orientable, such as for $\KU$ or $\KU\llpar q\rrpar$.

\begin{proof}
    We start with the commutative diagram of categories
    \[\begin{tikzcd}
        {\AVar_g(R)}\ar[d, "{[p^\infty]}"]   &    {\AVar_g(\tau_{\geq 0} R)}\ar[r]\ar[d, "{[p^\infty]}"]\ar[l, "\simeq", swap]   &   {\AVar_g(\pi_0 R)}\ar[d, "{[p^\infty]}"]   \\
        {\BTg(R)}   &   {\BTg(\tau_{\geq 0} R)}\ar[r]\ar[l, "\simeq", swap]   &   {\BTg(\pi_0 R)}
    \end{tikzcd}\]
    where all of the horizontal maps are induced by the functoriality of $\AVar_g(-)$ and $\BTg(-)$ in morphisms of $\E_\infty$-rings; see \cite[Pr.6.7.1]{ec1} for the vertical maps. The two left-most horizontal maps are equivalences by \cite[Rmks.1.5.3 \& 6.5.3]{ec1}, respectively; this is essentially a definition of the categories $\AVar_g(R)$ and $\BTg(R)$ when $R$ is not connective. The right-hand square is Cartesian by \cite[Pr.7.4.1]{ec1}, which concludes the proof.
\end{proof}

\begin{proof}[Proof of \Cref{orientedspectralserretate}]
    Consider the commutative diagram of categories
    \[\begin{tikzcd}
        {\Ell^\ori(R)}\ar[r]\ar[d, "{[p^\infty]}"]  &   {\Ell(R)}\ar[r]\ar[d, "{[p^\infty]}"]   &   {\Ell(\pi_0 R)}\ar[d, "{[p^\infty]}"]   \\
        {\BTtwoor(R)}\ar[r]  &   {\BTtwo(R)}\ar[r]   &   {\BTtwo(\pi_0 R)}
    \end{tikzcd}\]
    given by the obvious forgetful maps. By \Cref{spectralserretate}, the right square is Cartesian. Next, notice that for a fixed elliptic curve $E$ over $R$, the horizontal fibres of the left square above induce a map of spaces
    \[\Ori(E) \to \Ori(E[p^\infty])\]
    parametrising orientations of the elliptic curve $E$ and orientations of its associated $p$-divisible group $E[p^\infty]$. By \cite[Rmk.7.2.2]{ec2} and \cite[Pr.2.2.1]{ec3}, this map is an equivalence of spaces, hence the left square is Cartesian, which implies the result.
\end{proof}

The next key ingredient is the following corepresentability of $\TMF$.

\begin{theorem}[{\cite[Cor.4.6.2.11]{reconstruction}}]\label{corepresentabilityoftmfinnicecases}
    Let $\A$ be a complex-periodic $\E_\infty$-ring such that for some fixed finite integer $N$, the underlying classical Quillen formal group of $\A$ has height $\leq N$ for all primes $p$. Then the global sections functor induces a natural equivalence of spaces
    \[\Ell^\ori(\A)^\simeq = \Map_{\Stk}(\Spec \A,\M^\ori_\Ell) \xrightarrow{\simeq} \Map_{\CAlg}(\TMF,\A) .\]
\end{theorem}

We are now ready to prove \Cref{liftoftate}, which immediately implies \Cref{existenceoftatecurve}. To do this, we will use \Cref{orientedspectralserretate} to lift $\T_\km$ and the classical Tate curve to an elliptic curve over the $p$-completion of $\KU\llpar q \rrpar$ for each prime $p$. Then we will use \Cref{corepresentabilityoftmfinnicecases} and an arithmetic fracture square to glue these together with some additional rational data.

\begin{proof}[Proof of \Cref{liftoftate}]
    Let us write $\b=\KU\llpar q\rrpar$. For each fixed prime $p$, the $\P$-divisible group $\T_\km$ over $\b$ yields a $p$-divisible group $\widehat{\T}_\km$ over the $p$-completion $\widehat{\b}$ of $\b$. By \Cref{spectralliftofkm}, we see that the underlying classical $p$-divisible group of $\widehat{\T}_\km$ over $\pi_0 \widehat{\b} = \Z\llpar q\rrpar^\wedge_p$ is equivalent to the base-change of the Katz--Mazur $p$-divisible group $\T_\km^\heartsuit$ over $\pi_0 \widehat{\b}$. Katz--Mazur further identify this $p$-divisible group with the $p$-divisible group associated with the base-change of the classical Tate curve $\T^\heartsuit$ over $\pi_0 \widehat{\b}$; see \cite[(8.8)]{km}. By \Cref{orientedspectralserretate}, we have a Cartesian diagram of categories
    \[\begin{tikzcd}
        {\Ell^\ori(\widehat{\b})}\ar[r]\ar[d, "{[p^\infty]}"]   &   {\Ell(\pi_0 \widehat{\b})}\ar[d, "{[p^\infty]}"]   \\
        {\BTtwoor(\widehat{\b})}\ar[r]   &   {\BTtwo(\pi_0 \widehat{\b}).}
    \end{tikzcd}\]
    The oriented $p$-divisible group $\widehat{\T}_\km$ over $\widehat{\b}$ defines an object in the lower-left corner, the classical Tate curve $\widehat{\T}^\heartsuit$ over $\pi_0 \widehat{\b}$ defines an object in the upper-right corner, and the isomorphism identifying the classical $p$-divisible associated with $\widehat{\T}_\km$ over $\Z\llpar q\rrpar^\wedge_p$ of \Cref{spectralliftofkm} yields an isomorphism between their images in the lower-right corner---Katz--Mazur even state this isomorphism is unique; see property T.4 in \cite[(8.8)]{km}. This pullback of categories above then yields an oriented elliptic curve $\widehat{\T}$ over $\widehat{\b}$ for each prime $p$.\\

    Next, notice that $\b_\Q$, the rationalisation of $\b$, has homotopy groups
    \[\pi_\ast \b_\Q \simeq \Z\llpar q \rrpar[u^\pm] \otimes \Q,\qquad |u| = 2.\]
    We define an elliptic curve $\T_\Q$ over $\b_\Q$ as the base-change of the classical Tate curve $\T^\heartsuit$ over $\Z\llpar q \rrpar$ to this $\E_\infty$-$\Q$-algebra $\b_\Q$. This comes with a canonical orientation too as we are working rationally. The associated map of rational $\E_\infty$-rings $\TMF_\Q \to \b_\Q$ is precisely the rationalisation of (\ref{qexpmap}) on homotopy groups. The formal group associated with $\T_\Q$ is $\widehat{\G}_m$ as this is true for the classical Tate curve. 
    To glue each $\widehat{\T}$ with the rational Tate curve $\T_\Q$, we will use an arithmetic fracture square, so the Cartesian diagram of $\E_\infty$-rings
    \[\begin{tikzcd}
        {\b}\ar[r]\ar[d]  &   {\prod_p \widehat{\b}}\ar[d]  \\
        {\b_\Q}\ar[r]  &   {\left(\prod_p \widehat{\b}\right)_\Q.}
    \end{tikzcd}\]
    As the functor $\Ell^\ori(-)^\simeq\colon \CAlg\to \Spc$ is corepresented by $\TMF$ when restricted to complex periodic $\E_\infty$-rings by \Cref{corepresentabilityoftmfinnicecases}, we have equivalences of spaces
    \[\Ell^\ori(\b)^\simeq \simeq \Map_{\CAlg}(\TMF, \b) \simeq \Map_{\CAlg}(\TMF, \b_\Q)\underset{\Map_{\CAlg}(\TMF, (\prod_p \widehat{\b})_\Q)}{\times} \Map_{\CAlg}(\TMF, \prod_p \widehat{\b})\]
    \[\simeq \Ell^\ori(\b_\Q)^\simeq\underset{\Ell^\ori((\prod_p \widehat{\b})_\Q)^\simeq}{\times} \Ell^\ori(\prod_p \widehat{\b})^\simeq.\]
    We already have the desired oriented elliptic curves $\widehat{\T}$ over each $\widehat{\b}$ and $\T_\Q$ over $\b_\Q$, so to obtain $\T$ over $\b$, it suffices to produce an isomorphism between oriented elliptic curves over $(\prod_p \widehat{\b})_\Q$. To do this, consider the equivalences of spaces
    \[\Ell^\ori\left((\prod_p \widehat{\b})_\Q\right)^\simeq \simeq \Map_{\CAlg}\left(\TMF, (\prod_p \widehat{\b})_\Q\right) \simeq \Map_{\CAlg_\Q}\left(\TMF_\Q, (\prod_p \widehat{\b})_\Q\right)\]
    \begin{equation}\label{seconddecomposition}\subseteq \Map_{\CAlg_\Q}\left(\tmf_\Q, (\prod_p \widehat{\b})_\Q)\right) \simeq \Omega^{\infty+8}\left((\prod_p \widehat{\b})_\Q\right) \times \Omega^{\infty+12}\left((\prod_p \widehat{\b})_\Q\right),\end{equation}
    where we use that $\tmf_\Q$ is the free $\E_\infty$-$\Q$-algebra on two generators $c_4$ and $c_6$ in degrees $8$ and $12$, respectively, that $\TMF_\Q \simeq \tmf_\Q[\Delta^{-1}]$ for $\Delta = (c_4^3-c_6^2)/1728$, and that the containment above is a monomorphism onto those path components where the image of $\Delta$ is invertible. Both the base-change of $\widehat{\T}$ and $\T_\Q$ to $(\prod_p \widehat{\b})_\Q$ are chosen such that the associated map from $\tmf_\Q$ sends $c_4$ and $c_6$ to their $q$-expansion maps. From the equivalences above, this shows that these oriented elliptic curves are equivalent over $(\prod_p \widehat{\b})_\Q$. The fact that $(\prod_p \widehat{\b})_\Q$ has homotopy groups concentrated in even degrees shows that the choice of equivalence is unique up to homotopy; even up to $3$-homotopy. This yields the oriented elliptic curve $\T$ over $\b$.\\
    
    To compute the formal group associated with $\T$, notice that $\b$ is an $\E_\infty$-$\KU$-algebra, as it then follows that the formal group associated with $\T$ is $\widehat{\G}_m$. Indeed, by \cite[Pr.4.3.23]{ec2}, an orientation on $\T$ is an equivalence between its associated formal group and the Quillen formal group over $\b$. As the Quillen formal group of a complex-periodic $\E_\infty$-ring satisfies base-change, we obtain the desired result of the Quillen formal group $\widehat{\G}_m$ over $\KU$. Alternatively, one can use the fact that the space of oriented formal groups over a complex-periodic $\E_\infty$-ring is contractible.\\
    
    To identify the $\P$-divisible group associated with $\T$ with the base-change of $\T_\km$ to $\b$, we compute its $N$-torsion points $\T[N]$ for each integer $N\geq 2$. As $\T$ is an elliptic curve over $\b$, each $\T[N]$ is a finite flat scheme over $\b$, so of the form $\Spec \b_N$ for some finite flat $\E_\infty$-$\b$-algebra $\b_N$. By an arithmetic fracture square, it suffices then to identify $\b_N$ with the $N$-torsion points of $\T_\km$ over $\widehat{\b}$ for each prime $p$, over $\b_\Q$, and finally over $(\prod_p \widehat{\b})_\Q$. This then follows from the fact that $N$-torsion of an elliptic curve is naturally equivalent to the $N$-torsion of its associated $\P$-divisible group together with the construction of $\T$ above; this fact is discussed in \cite[Pr.4.3]{elltempcomp}. Finally, the uniqueness of $\T$ also follows from the proof above, in particular, from the computation
    \[\pi_1 \Ell^\ori\left((\prod_p \widehat{\b})_\Q\right) \simeq \pi_1 \Map_{\CAlg_\Q}\left(\TMF_\Q, (\prod_p \widehat{\b})_\Q\right)=0\]
    by (\ref{seconddecomposition}) and the evenness of $(\prod_p \widehat{\b})_\Q$.
\end{proof}



\section{Compactifying the moduli stack of oriented elliptic curves}\label{sec:compact}

Given the derived Tate curve $\T$ over $\Spec \KU\llpar q\rrpar$ of \Cref{liftoftate}, we can now follow the construction of a compactification of $\M_\Ell^\ori$ outlined in \cite[\textsection4.3]{lurieecsurvey}. First, we need to account for some automorphisms of $\T$.\\

There is a $C_2$-action on $\Spec \KU\llpar q \rrpar$ given by complex conjugation. This is a wonderful description of this action, as it comes from a geometric description and is useful for computations, however, we would like to reinterpret this $C_2$-action to more easily work in spectral algebraic geometry. \\

First, notice that $\KU\llpar q \rrpar$ sits in the pushout of $\E_\infty$-rings
\[\begin{tikzcd}
    {\Sph}\ar[r]\ar[d]              &   {\KU}\ar[d] \\
    {\Sph\llpar q \rrpar}\ar[r]     &   {\KU\llpar q \rrpar.}
\end{tikzcd}\]
In the language of \cite[Df.4.3.13]{ec2}, we view $\KU$ as the orientation classifier of $\widehat{\G}_m$ over $\Sph$, meaning that the universal orientation $e$ of $\KU$ induces an equivalence of spaces
\[\Map_{\CAlg}(\KU,\KU) \xrightarrow{\simeq} \Ori(\widehat{\G}_{m,\KU});\]
see \cite[Th.6.5.1]{ec2}. By \cite[Rmk.6.4.2]{ec2}, orientation classifiers are preserved by base-change, meaning that the pushout of $\E_\infty$-rings above shows that $\KU\llpar q \rrpar$ is the orientation classifier of $\widehat{\G}_m$ over $\Sph\llpar q \rrpar$. To reiterate, the universal orientation of $\widehat{\G}_m$ over $\KU\llpar q\rrpar$ induces the equivalence of spaces
\begin{equation}\label{eq:upofkuq}\Map_{\CAlg_{\Sph\llpar q \rrpar}}(\KU\llpar q \rrpar,\KU\llpar q \rrpar) \xrightarrow{\simeq} \Ori(\widehat{\G}_{m,\KU\llpar q \rrpar}).\end{equation}
As $\T$ is oriented over $\KU\llpar q\rrpar$, there is a preferred identification of formal groups $\widehat{\T} \simeq \widehat{\G}_m$ over $\KU\llpar q\rrpar$. We then equip the right-hand side of (\ref{eq:upofkuq}) with a $C_2$-action on the right-hand side via the involution $[-1]$ of $\T$, yielding a involution of $\E_\infty$-$\Sph\llpar q \rrpar$-algebras
\begin{equation}\label{eq:ctwoaction} c\colon \KU\llpar q \rrpar \to \KU\llpar q \rrpar.\end{equation}
By linearity, $c(q)=q$. Consider the diagram of stacks
\[\begin{tikzcd}
    {\Spec \KU\llpar q \rrpar}\ar[dr, "{(\T,e)}", swap]\ar[rr, "{c}"] &&  {\Spec \KU\llpar q \rrpar}\ar[dl, "{(\T,e)}"]   \\
        &   {\M_\Ell^\ori.}  &
\end{tikzcd}\]
By definition, the composite $(\T,e)\circ c$ is represented by the pair $(c^\ast \T=\T, [-1](e))$. Inside $\M_\Ell^\ori(\Spec\KU\llpar q \rrpar)$, the pairs $(\T, e)$ and $(\T,[-1]e)$ are equivalent via the isomorphism of elliptic curves $[-1]\colon \T\to \T$. In other words, this isomorphism of $\T$ induces a filling for the diagram above. This naturality of $[-1]$---the fact it naturally commutes with all isomorphisms of oriented elliptic curves---implies that $\T\colon \Spec \KU\llpar q \rrpar \to \M_\Ell^\ori$ descends to the desired map
\begin{equation}\label{mapfactoringthroughquotient}
\T \colon \Spec \KU\llpar q\rrpar/ C_2 \to \M_\Ell^\ori
\end{equation}
through the quotient of $\Spec \KU\llpar q \rrpar$ by this $C_2$-action. In \cite[\textsection A]{geometricnorms}, we will show how $\Spec \KU\llpar q\rrpar / C_2$ can be interpreted as the \emph{moduli stack of forms of the oriented smooth Tate curve} $\M_\Tate^\ori$; this will not play a role here.

\begin{remark}[Relation to complex conjugation]\label{complexconjguation}
This $C_2$-action $[-1]$ on $\KU\llpar q\rrpar$ is given by the stable Adams operation $\psi^{-1}$ on $\KU\llpar q\rrpar$. To be precise, as $\Z \to \Z\llpar q\rrpar$ is a flat map of ring, the cohomology theory $\KU\llpar q\rrpar^\ast(-)$ is computed by basechanging $\KU^*$ to $\Z\llpar q\rrpar$:
\begin{equation}\label{tatektheorycohomologyhteory}\KU\llpar q\rrpar^\ast(-) \simeq \KU^\ast(-) \otimes_{\Z} \Z\llpar q\rrpar \qquad |q|=0.\end{equation}
From this construction, Tate $K$-theory obtains the Adams operation $\psi^{-1}$ from $\KU$ inherited from the complex conjugation operation on complex vector bundles. To compare these two $C_2$-actions as natural ring operations on cohomology theories, we can use the splitting principle and additivity to reduce ourselves to computing their action on the universal line bundle $\ga$ over $\B U(1)$. The computation for $\psi^{-1}$ this is well-known, and for $[-1]$ this comes from the same arguments as made in the proof of \cite[Pr.6.21]{luriestheorem}; an integral version of this will be discussed in \cite[\textsection A]{geometricnorms}. Due to the lack of phantom maps between Landweber exact cohomology theories, see \cite[Cor.2.15]{hoveysticklandMoravalocal}, this induces a homotopy between maps of spectra. These $C_2$-actions also agree as endomorphisms of the $\E_\infty$-ring $\KU\llpar q\rrpar$, which can be proven using Goerss--Hopkins obstruction theory and an arithmetic fracture square; we leave those details to the interested reader.
\end{remark}

\begin{remark}[Global sections]\label{mapfromtmftoko}
    Further refining \Cref{complexconjguation}, we can identify the global sections of $\Spec \KU \llpar q \rrpar /C_2$ with $\KO\llpar q \rrpar$, just as the global sections of $\Spec \KU /C_2$ can be identified with $\KO$. In particular, we now have morphisms of $\E_\infty$-rings
    \[\TMF \to \KO\llpar q \rrpar \to \KU\llpar q \rrpar\]
    by taking global sections of (\ref{mapfactoringthroughquotient}) and the quotient $\Spec \KU\llpar q \rrpar \to \Spec \KU \llpar q \rrpar / C_2$.
\end{remark}

\begin{remark}[Stable Adams operations]\label{rmk:adamsoperations}
    Let us work over $\Sph[\tfrac{1}{k}]$ for some integer $k$. Then the base-change of (\ref{mapfactoringthroughquotient}) yields a map of stacks
    \[\T \colon \Spec \KU[\tfrac{1}{k}]\llpar q\rrpar/ C_2 \to \M_\Ell^\ori \times \Spec \Sph[\tfrac{1}{k}].\]
    Using the obvious variation of (\ref{eq:upofkuq}), one can construct stable Adams operations $\psi^{k}$ on $\KU[\tfrac{1}{k}]\llpar q \rrpar$ via the $k$-fold multiplication map $[k]\colon \T \to \T$, which now induces an automorphism of $\widehat{\G}_m$ over $\KU[\tfrac{1}{k}]\llpar q \rrpar$. By construction, or rather by their characterisation from (\ref{eq:upofkuq}), these Adams operations are just as coherent as the multiplication maps; this is discussed in \cite[Rmk.2.16]{heckeontmf}. In particular, we have natural homotopies between the automorphisms of $\E_\infty$-rings $\psi^{k}\circ \psi^\ell \simeq \psi^{k\ell}$ for integers $k,\ell$. Notice that as these maps $\psi^k$ are induced by $[k]\colon \T\to \T$ which is \textbf{not} an equivalence of elliptic curves for $k\neq \pm1$, the automorphisms $\psi^k$ are not $\TMF$-linear.
\end{remark}

To construct $\overline{\M}_\Ell^\ori$, we will use the following fact from spectral algebraic geometry.

\begin{lemma}\label{pushoutexists}
    The pushout of the cospan
    \begin{equation}\label{cospanforpo}\M_\Ell^\ori \xleftarrow{\T} \Spec \KU\llpar q \rrpar / C_2 \to \Spec \KU\llbracket q \rrbracket / C_2\end{equation}
    of $1$-localic stacks exists.
\end{lemma}

\begin{proof}
    First, we take the pushout in the category of $1$-localic spectrally ringed $\infty$-topoi
    \begin{equation}\label{zerothpushouttopi}
        \begin{tikzcd}
            {\Spec \KU\llpar q \rrpar / C_2}\ar[d, "{\T}"]\ar[r]    &   {\Spec \KU \llbracket q \rrbracket / C_2}\ar[d, "g"]   \\
            {\M_\Ell^\ori}\ar[r, "f"]    &   {(X,\O_X),}
        \end{tikzcd}
    \end{equation}
    see \cite[Df.1.4.2.1]{sag}, which makes sense as $\M_\Ell^\ori$ is $1$-localic by \cite[Th.2.4.1]{ec1} and \cite[Rmk.7.3.2]{ec2}. In particular, the underlying $1$-localic $\infty$-topos $X$ can be computed as $\Shv(X^\heartsuit)$. As the truncation functor from $\infty$-topoi to $1$-topoi is a left adjoint, the diagram of classical topoi
    \begin{equation}\label{firstpushoutpotoi}\begin{tikzcd}
        {\Spec \Z\llpar q \rrpar / C_2}\ar[r, "{h^\heartsuit}"]\ar[d] &   {\Spec \Z\llbracket q \rrbracket / C_2}\ar[d]   \\
        {\M_\Ell}\ar[r]  &   {(X^\heartsuit, \pi_0 \O_X)}
    \end{tikzcd}\end{equation}
    is a pushout, which recognises $X^\heartsuit$ as the classical $1$-topos underlying the étale site on the compactification of the moduli stack of elliptic curves. Moreover, the pushout (\ref{zerothpushouttopi}) yields a Meyer--Vietoris sequence of sheaves of spectra on $X$
    \begin{equation}\label{keymvfibresequence}\O_X \to f_\ast \O_{\M_\Ell^\ori} \times g_\ast \O_{\Spec \KU\llbracket q \rrbracket / C_2} \to f_\ast \T_\ast \O_{\Spec \KU\llpar q \rrpar / C_2}.\end{equation}
    One of the defining features of $\M_\Ell^\ori$, see \cite[Rmk.7.3.2]{ec2} for example, and also $\Spec \KU\llbracket q \rrbracket / C_2$ and $\Spec \KU\llpar q \rrpar / C_2$ is the fact that these complex-oriented stacks have even structure sheaf, so each of these stacks admits an étale cover by affines $\Spec A$ where the $\E_\infty$-ring $A$ has homotopy groups concentrated in even degrees. In particular, the odd homotopy group sheaves associated with (\ref{keymvfibresequence}) vanish, giving us the short exact sequence of sheaves of abelian groups on $X^\heartsuit$
    \begin{equation}\label{eq:mvsequenceonhomotopygroups}0 \to \pi_{2n} \O_{X} \to \pi_{2n} f_\ast \O_{\M_\Ell^\ori} \times \pi_{2n} g_\ast \O_{\Spec \KU \llbracket q \rrbracket / C_2} \to \pi_{2n} f_\ast T_\ast \O_{\Spec \KU\llpar q \rrpar / C_2} \to 0.\end{equation}
    Setting $n=0$, this exact sequence witnesses $\pi_0 \O_X = \O_{X^\heartsuit}$ as the classical structure sheaf for the compactification of the classical moduli stack of elliptic curves, therefore identifying $(X^\heartsuit, \pi_0 \O_X)$ with the classical compactification of the moduli stack of elliptic curves itself. 

    By \cite[Th.1.4.8.1]{sag}, to conclude that $X$ itself is a stack, ie, that it is Deligne--Mumford, it suffices to show that the homotopy sheaves $\pi_n \O_X$ are quasi-coherent sheaves over $(X^\heartsuit, \pi_0 \O_X)$ and that $\O_X$ itself is hypercomplete. The first follows from (\ref{eq:mvsequenceonhomotopygroups}), as quasi-coherent sheaves are closed under taking kernels, and the second from (\ref{keymvfibresequence}), as hypercomplete sheaves are closed under taking fibres.
\end{proof}

\begin{mydef}
Let us write $\overline{\M}_\Ell^\ori$ for the stack of \Cref{pushoutexists}, in other words, defined by the pushout
    \begin{equation}\label{thediagram}\begin{tikzcd}
    {\Spec \KU\llpar q\rrpar/C_2}\ar[d,"{\T}"]\ar[r]    &   {\Spec \KU\llbracket q\rrbracket/C_2}\ar[d] \\
    {\M_{\Ell}^\ori}\ar[r]  &   {\overline{\M}_{\Ell}^\ori.}
\end{tikzcd}\end{equation}
\end{mydef}

Some properties of $\overline{\M}_\Ell^\ori$ follow immediately from our considerations above.

\begin{prop}\label{basicproperties}
    The spectral Deligne--Mumford stack $\overline{\M}_\Ell^\ori$ is proper, almost of finite presentation, and $1$-localic with underlying classical Deligne--Mumford stack $\overline{\M}_\Ell$ the classical compactification of the moduli stack of elliptic curves of \cite{delignemumford}. Writing $\O_{\overline{\M}_\Ell^\ori}$ for the structure sheaf of $\overline{\M}_\Ell^\ori$, then for each integer $n$ there are isomorphisms of sheaves of abelian groups over $\overline{\M}_\Ell$
    \[\pi_{2n}\O_{\overline{\M}_\Ell^\ori} \simeq \omega_{\overline{\M}_\Ell}^{\otimes n}, \qquad \pi_{2n+1}\O_{\overline{\M}_\Ell^\ori} = 0.\]
\end{prop}

\begin{proof}
    The identification of the underlying classical Deligne--Mumford stack associated with $\overline{\M}_\Ell^\ori$ and the fact that its $1$-localic were shown in \Cref{pushoutexists}; in particular, it was shown that (\ref{firstpushoutpotoi}) was a pushout of classical Deligne--Mumford $1$-topoi. A map of spectral Deligne--Mumford stacks $\x \to \y$ is proper is $\tau_{\geq 0}\x \to \tau_{\geq 0}\y$ is proper à la \cite[Df.5.1.2.1]{sag}. By \cite[Rmk.5.1.2.2]{sag}, our desired map is proper if and only if the associated map of underlying classical Deligne--Mumford stacks is proper, which is true classically. The fact that $\overline{\M}_\Ell^\ori$ is almost of finite presentation over $\Spec \Sph$ follows by pushing (\ref{keymvfibresequence}) down to the $\Sph$. The identificiations of homotopy groups come from (\ref{eq:mvsequenceonhomotopygroups}).
\end{proof}

The final thing we would like to check in this section is that the above definition of $\overline{\M}_\Ell^\ori$ agrees with the previous definition of Goerss--Hopkins--Miller which we will denote by $\overline{\M}_\Ell^\GHM$; see \cite[\textsection12]{tmfbook} by Behrens for this construction. We can now finish our proof of \Cref{comparison}.

\begin{theorem}\label{comparisonintext}
    There is an equivalence of spectral Deligne--Mumford stacks
    \[\overline{\M}_\Ell^\ori \simeq \overline{\M}_\Ell^\GHM.\]
	In particular, there is an equivalence of $\E_\infty$-rings $\Ga(\overline{\M}_\Ell^\ori) \simeq \Tmf$ and a map of $\E_\infty$-rings $\Tmf \to \KO\llbracket q \rrbracket$, which we call the \emph{holomorphic topological $q$-expansion map}.
\end{theorem}

This allows us to \emph{define} $\Tmf$ as the global sections of $\overline{\M}_\Ell^\ori$ and a canonical holomorphic topological $q$-expansion map. Post-composing this with evaluation at $q=0$ gives the famous map of $\E_\infty$-rings $\Tmf \to \KO$ comparing topological modular forms and topological $K$-theory. It is crucial to use $\Tmf$ to define this direct comparison to $\KO$; this map does not factor through $\TMF$ since it sends $\Delta^{24}$ to zero, as $\Delta$ is a cusp form.

\begin{proof}
    Both stacks have the same underlying classical stack. As both stacks are $1$-localic, it remains to compare their structure sheaves. By the main theorem of \cite{uniqueotop}, it suffices to check that the structure sheaf of $\overline{\M}_\Ell^\ori$, which we write as $\O$, is naturally a sheaf of Landweber exact elliptic cohomology theories on these affines. Indeed, the statement of \cite[Th.1.3]{uniqueotop} is that all such sheaves of $\E_\infty$-rings with this extra data and properties are equivalent to $\O^\top$, the latter being the structure sheaf of $\E_\infty$-rings for $\overline{\M}_\Ell^\GHM$.\\
    
    To check that for each étale morphism $E\colon \Spec R \to \overline{\M}_\Ell$, the $\E_\infty$-ring $\O(\Spec R) = \A$ has the structure of an elliptic cohomology theory for $E/\Spec R$ natural in $R$ comes down to checking the four conditions of \cite[Df.1.1]{uniqueotop}. Parts (1), (2), and (3) are clear, as we have identified the homotopy groups of $\A$ in \Cref{basicproperties}. For part (4), we refer to \cite[Lm.4.5]{hilllawson}, which states that pullbacks of elliptic cohomology theories are themselves elliptic cohomology theories, which by (\ref{keymvfibresequence}) gives us the desired result; note that in \cite[Lm.4.5]{hilllawson}, the authors forgot to ask for the cospan to induce a jointly surjective map on $\pi_0$, however, this is clearly true in our case by (\ref{eq:mvsequenceonhomotopygroups}).
\end{proof}

An obvious consequence of the above comparison is that many of the now classical facts about $\overline{\M}_\Ell^\GHM$ also follow for $\overline{\M}_\Ell^\ori$. For example, one can compute the homotopy groups of $\Tmf$ using the \emph{descent spectral sequence}
\[E_2^{s,t} = H^s(\overline{\M}_\Ell, \omega^{\otimes t/2}) \implies \pi_{t-s} \Tmf\]
as done in \cite{smfcomputation}. The affineness results of Mathew--Meier \cite[Th.7.2(2)]{akhilandlennart}, also see \cite[Ex.4.2.3.5 \& Cor.4.4.0.3]{reconstruction}, imply that the global sections functor
\[\Ga\colon \QCoh(\overline{\M}_\Ell^\ori) \xrightarrow{\simeq} \Mod_{\Tmf}\]
is an equivalence of symmetric monoidal categories. This also allows us to identify the map of $\E_\infty$-rings $\Tmf \to \TMF$ associated with the open immersion of stacks $\M^\ori_\Ell \to \overline{\M}_\Ell^\ori$ with the localisation map $\Tmf \to \Tmf[\Delta^{-24}]$, for a certain element $\Delta^{24}\in \pi_{576}\Tmf$.


\section{Equivariant and global Tate \texorpdfstring{$K$}{K}-theory}\label{sec:globalthings}

Having used equivariant homotopy theory to produce the derived Tate curve $\T$ and the derived Katz--Mazur $\mathbf{P}$-divisible group $\T_\km$, we now switch focus and use these objects to construct various genuine global cohomology theories. \\

In the first subsection, we use the derived Tate curve to construct the global smooth Tate $K$-theories $\gKO\llpar q\rrpar$ and $\gKU\llpar q\rrpar$, as well as a morphism of global $\E_\infty$-rings $\gTMF \to \gKO\llpar q \rrpar$, thus proving \Cref{globalmap}. In the second subsection, we use the derived Katz--Mazur $\mathbf{P}$-divisible group to define \emph{global quasi-elliptic cohomology} $\gQEll$ and prove some of its basic properties summarised in \Cref{cor:qellstuffintro}. In short, this first subsection studies $\T$ and its associated equivariant cohomology theory for compact Lie groups, and the second subsection studies $\T_\km$ and its equivariant cohomology theory for finite groups.\\

As a final application, we consider the \emph{$C_2$-equivariant global smooth Tate $K$-theory} and its relationship with the $C_2$-equivariant global spectrum $\gTMF_1(3)$ of \cite{levelonethree} and \cite{lennartandhill}, proving \Cref{cor:realstuffintro}.


\subsection{Global smooth Tate \texorpdfstring{$K$}{K}-theory}\label{ssec:globaltate}

Recall that $\Sp_{\gl}$ denotes the category of \emph{global spectra} in the sense of Schwede. For every compact Lie group $G$ there is a forgetful functor $\mathrm{res}_G\colon\Sp_{\gl}\to \Sp_G$. Given a global spectrum $\gX$ we will typically write $\gX_G$ for the $G$-spectrum associated with $\gX$. In particular, for $G=e$ the trivial group, we obtain a spectrum $\gX_e$, the \emph{underlying spectrum} of $\gX$. \\

There is also a \emph{global suspension spectrum} functor $\Sigma_+^\infty \colon \Spc_\gl \to \Sp_\gl$. Using this, any global spectrum induces a cohomology theory on global spaces via the functor $\Map_{\Sp_{\gl}}(\Sigma_+^\infty(-),\gX).$
Such equivariant cohomology theories come equipped with significantly more structure than those of \cite{ec3} or \cite{davidandlennart}, as the former have dimension shifting transfers.\\

Given a stack $\s$ and $\G$ an oriented geometric abelian $\s$-stack (\Cref{defofgeometricstacks}), we obtain an equivariant cohomology theory on global spaces via Gepner--Meier's \emph{equivariant elliptic cohomology} construction of \cite[Con.6.1]{davidandlennart}. This yields a colimit preserving functor
\[
\Ell_{\G/\s}\colon \Spc_\gl \to \Shv(\s).
\]
The quasi-coherent sheaves $\Ell_{\G/\s}(\B G)$ are representable by stacks over $\s$ if $G$ is an abelian compact Lie group. Applying the global sections functor $\Gamma\colon \Shv(\s)\to \Sp$ gives a global cohomology theory.\\

We can now state our simplified version of \cite[Th.D]{gepner2024global2ringsgenuinerefinements}.

\begin{theorem}\label{functorialityoftwoglobal}
    Let $\s$ be a stack and $\G$ be an oriented elliptic curve over $\s$. Then there is a functor
    \[\Ga_{\gl} \colon \Stk_{/\s}^\op \to \CAlg(\Sp_\gl)\]
    such that for an $\s$-stack $f\colon \t \to \s$, the underlying $\E_\infty$-ring of $\Ga_\gl(\t)^{\Sigma_+^\infty \B H}$ is given by the global sections of $\Ell_{f^\ast\G/\t}^{\B H}$ for an \emph{abelian} compact Lie group $H$.
\end{theorem}

\begin{remark}\label{extend_from_ab}
The global spectra of the theorem above are constructed first as global spectra for the family of abelian compact Lie groups. We then apply the right adjoint to the restriction of families to obtain a global spectrum for the family of all compact Lie groups. See \cite[Ex.14.23]{gepner2024global2ringsgenuinerefinements} for more details. One should compare this with the global spectrum representing quasi--elliptic cohomology constructed in the next section, which has to be extended from \emph{finite} abelian groups.
\end{remark}


The functoriality of the theorem is immediate from the functoriality of the construction of $\QCoh(\Ell_{(-)})$ from \cite{davidandlennart} in base-change combined with the functoriality of \cite[Pr.12.8]{gepner2024global2ringsgenuinerefinements}.\\

Using \Cref{functorialityoftwoglobal}, Gepner--Pol and the second-named author define global versions of equivariant topological modular forms.

\begin{mydef}\label{globaltmfdefintion}
    The global $\E_\infty$-ring of \emph{global topological modular forms} $\gTMF$ is the global $\E_\infty$-ring given by applying $\Ga_\gl$ to the universal oriented elliptic curve $\e^\ori$ living over $\M_\Ell^\ori$.
\end{mydef}

We can also apply \Cref{functorialityoftwoglobal} to the oriented elliptic curve $\T$ from \Cref{liftoftate}, which yields global versions of equivariant Tate $K$-theory.

\begin{mydef}\label{globaltatedefinition}
    Define the global $\E_\infty$-rings $\gKO\llpar q \rrpar$ and $\gKU\llpar q \rrpar$ by applying the functor $\Ga_{\gl}$ of \Cref{functorialityoftwoglobal} to $\T$ over $\Spec \KU\llpar q \rrpar / C_2$ and $\Spec \KU\llpar q \rrpar$, respectively. We call these global $\E_\infty$-rings \emph{global smooth real Tate $K$-theory} and \emph{global smooth complex Tate $K$-theory}, respectively. By \Cref{functorialityoftwoglobal}, this comes with a map of global $\E_\infty$-rings $\gKO\llpar q \rrpar \to \gKU \llpar q\rrpar$ refining the map of $\E_\infty$-rings $\KO \llpar q \rrpar \to \KU\llpar q \rrpar$. 
\end{mydef}

We use the adjective \emph{smooth} above to distinguish between $\gKU\llpar q \rrpar$ and a potential $\gKU\llbracket q \rrbracket$ which we do not define here; over $\Z\llpar q \rrpar$ the classical Tate curve is a \emph{smooth} elliptic curve and over $\Z\llbracket q \rrbracket$ it has a nodal singularity. 

\begin{warn}
    While Tate $K$-theory $\KU\llpar q \rrpar$ can be expressed as a cohomology theory using (\ref{tatektheorycohomologyhteory}), there is \emph{not} an equivalence between $\gKU\llpar q \rrpar^\ast(X)$ and $\gKU^\ast(X) \otimes_{\Z} \Z\llpar q \rrpar$ for arbitrary global spaces $X$. Indeed, the latter is a \emph{globally flat} cohomology theory over $\gKU$, using the terminology of \cite[Df.3.1]{firstpreprint}. Since the Tate curve is not base-changed from $\G_m$, this is not the case for $\gKU\llpar q \rrpar$.   
\end{warn}

The following is a slight refinement of \Cref{globalmap}.

\begin{cor}
    There are maps of global $\E_\infty$-rings
    \begin{equation}\label{globalmaprefined}\gTMF \to \gKO\llpar q \rrpar \to \gKU\llpar q \rrpar\end{equation}
    whose underlying maps of $\E_\infty$-rings are those from \Cref{mapfromtmftoko}. In particular, the map of rings induced by taking $\pi_\ast^e(-)$ is the composition of the classical meromorphic $q$-expansion map (\ref{qexpmap}) through the edge map $\pi_\ast \TMF \to \MF$ in the descent spectral sequence.
\end{cor}

\begin{proof}
    This follows from \Cref{functorialityoftwoglobal} as the map of stack $f\colon \Spec \KU\llpar q \rrpar / C_2 \to \M_\Ell^\ori$ tautologically defines $\T$ as the pullback of the universal oriented elliptic curve along $f$. This also goes for the quotient map $\Spec \KU\llpar q \rrpar \to \Spec \KU \llpar q \rrpar / C_2$.
\end{proof}

As equivariant elliptic cohomology theories, both $\gKO\llpar q \rrpar$ and $\gKU\llpar q \rrpar$ inherit many of the good properties discussed in \cite{davidandlennart}. For example, for any compact Lie group $G$ of the form $F\times H$ where $F$ is a finite group and $H$ is a torus, then the $G$-fixed points of these theories all define dualisable $\KO\llpar q \rrpar$- and $\KU\llpar q \rrpar$-modules, respectively; see \cite[Th.C]{elltempcomp}. This does not hold for $\gKO$ or $\gKU$ when $G$ is a torus, for example. By \cite[Th.4.9.2]{ec3} and \cite[Th.A \& B]{elltempcomp}, $\gKU$ also satisfies the Atiyah--Segal completion theorem for any finite group $G$.


\subsection{Global quasi-elliptic cohomology}\label{ssec:qell}
Lurie's \emph{tempered cohomology} construction of \cite[Con.4.0.3]{ec3} applied to a $\P$-divisible group $\G$ yields a limit preserving functor
\[\Ga\Temp_{\G/\s} \colon \Spc_{\ab\fin}^\op \to \CAlg_{\Ga(\s)}.\]
The affine case $\s = \Spec \A$ was originally constructed in \cite{ec3}, where it is written as $\A_{\G}$; the extension to general stacks $\s$ appears in \cite[\textsection2]{elltempcomp} and \cite[\textsection4.6]{tempered}.\\

Recall that global spectra can be defined with respect to any family of compact Lie groups. As a special case, we may consider $\Sp_{\fingl}$, the category of \emph{finite} global spectra. Intuitively, a fin-global spectrum $\gX$ is a collection of $G$-spectra $\gX_G$ for each \emph{finite} group $G$ together with the usual compatibility. More formally, we define $\Sp_{\fingl}$ analogously to $\Sp_{\gl}$ as a partially lax limit, but now the indexing diagram is restricted to only include finite groups.\\

We now state our simplified version of \cite[Thm.A]{tempered}; the similar statement \cite[Thm.E]{gepner2024global2ringsgenuinerefinements} also suffices for many of our purposes here.

\begin{theorem}\label{functorialityoftwoglobalfin}
Given an oriented $\P$-divisible group $\G$ over an affine stack $\s = \Spec \A$, then there is a functor
\[\Ga_{\temp} \colon \CAlg_{\A/} \to \CAlg(\Sp_{\fingl})\]
such that for every $\E_\infty$-$\A$-algebra $f\colon \A \to \b$, the underlying $\E_\infty$-ring of $\Ga_{\temp}(\b)^{\Sigma_+^\infty \B H}$ is given by $\Ga \Temp_{f^*\G[\P^\infty]/\Spec(\b)}^{\B H}$ for a \emph{finite abelian} group $H$.
\end{theorem}

\begin{remark}
Compared to \Cref{functorialityoftwoglobal}, in this case, we first construct a global spectrum for the family of \emph{finite} abelian groups, which we then extend to all finite groups as in \Cref{extend_from_ab}. While one could extend the fin-global spectra $\Gamma_{\temp}(\b)$ to all compact Lie groups, we are not sure if the values of this extension on compact Lie groups are very natural objects to consider and so we prefer to restrict tempered cohomology to finite isotropy. Let us note however that $\gKU$ is even extended from \emph{finite cyclic groups}, so in this case, the result is a posteriori quite reasonable.
\end{remark}

\begin{warn}\label{finiteandnonfiniteremark}
By definition any oriented geometric abelian $\s$-stack $\G$ restricts to a $\mathbf{P}$-divisible group $\G[\P^\infty]$. In particular, we should compare the fin-global ring $\Ga_{\temp}\G[\P^\infty]$ with the restriction of  $\Ga_{\gl}\G$ to a fin-global ring. The main theorems of \cite{elltempcomp} state that the underlying fin-globally equivariant cohomology theories $\Ell_{\G/\s}$ and $\Ga\Temp_{\G[\P^\infty]/\s}$ agree. However, this does \textbf{not} immediately imply that the global refinements of \Cref{functorialityoftwoglobal} necessary agree. Indeed, to phrase this in the language of \cite{gepner2024global2ringsgenuinerefinements}, it is easy to construct a comparison between na\"{i}ve global $2$-rings using the results of \cite{elltempcomp}, but the two associated genuine global $2$-rings are not genuine with respect to the same global family. This means that the conjecture of \cite[Con.0.1]{elltempcomp} remains open, although we expect it to be true.
\end{warn}

One could also use the full power of \cite[Th.A]{tempered} in \Cref{functorialityoftwoglobalfin}, ie, the functor could takes values in $\Sp^\gl_\pi$, the category of \emph{$\pi$-ambidextrous global spectra}, where one has extra transfer maps available; see \cite[Ex.5.1.6-8]{tempered}. As this extension does not provide any significant advantages here, and the forgetful functor $\Sp^\gl_\pi \to \Sp_{\fingl}$ is conservative, we forgo that discussion here.\\

In this subsection, we introduce and summarize some properties of the fin-global $\E_\infty$-ring associated with the lift of the Katz--Mazur $\P$-divisible group $\T_\km$ to $\KU[q^\pm]$ of \Cref{spectralliftofkm}. The homotopy groups of this global spectrum have slightly nicer finiteness properties than $\gKU\llpar q\rrpar$ by virtue of being defined over $\KU[q^\pm]$.

\begin{mydef}
Let $\gQEll$ be the global $\E_\infty$-ring obtained by applying $\Ga_{\temp}$ to $\T_\km$ (over $\KU[q^\pm])$ and call it \emph{global quasi-elliptic cohomology}.
\end{mydef}

Equivariant quasi-elliptic cohomology was implicitly introduced by Ganter \cite{ganterstringy}, explicitly defined by Rezk in \cite{rezkquasielliptic}, and further explored by Huan in \cite{huanquasielliptic}. We call the formulation of Huan \emph{GRH-quasi-elliptic cohomology}. By \Cref{rmk:compwithgrh}, the rings $\pi_\ast^G \gQEll$ agree with the coefficients of GRH-quasi-elliptic cohomology for finite abelian groups $G$.\\

We would like to prove a structured version of one of the characterising properties of $\gQEll$: its comparison to global smooth Tate $K$-theory. To make this precise, we will need a finite-global version of smooth Tate $K$-theory. 

\begin{mydef}
Similar to \Cref{globaltmfdefintion,globaltatedefinition}, we can apply $\Ga_{\temp}$ to the oriented $\P$-divisible group $\T[\P^\infty]$ of $\T$ over $\Spec \KU\llpar q \rrpar$, which produces a fin-global $\E_\infty$-ring $\gKU_\fin \llpar q \rrpar$. As mentioned in \Cref{finiteandnonfiniteremark} we expect, but do not know, that this agrees with $\gKU\llpar q \rrpar$ after the latter is restricted to finite groups.
\end{mydef}

The functor $(-)_e\colon \Sp_{\fingl} \to \Sp$ is strong symmetric monoidal and also has a strong symmetric monoidal left adjoint $L$; see \cite[Th.4.5.1(3)]{s}. As $\gQEll$ has underlying $\E_\infty$-ring $\KU[q^\pm]$, the counit of this adjunction gives us a morphism of fin-global $\E_\infty$-rings $L \KU[q^\pm] \to \gQEll$. Applying this to $\gKU_\fin\llpar q \rrpar$, we also obtain a map $L \KU\llpar q \rrpar \to \gKU_\fin\llpar q \rrpar$. Using these maps, we can compare $\gQEll$ with $\gKU_\fin\llpar q\rrpar$.

\begin{theorem}\label{basechangeprop}
    The natural map of fin-global $\E_\infty$-rings
    \begin{equation}\label{bcofqell}\gQEll \otimes_{L \KU[q^\pm]} L \KU\llpar q \rrpar \xrightarrow{\simeq} \gKU_\fin\llpar q \rrpar\end{equation}
    is an equivalence. In particular, for any finite group $G$, the ring homomorphism
    \[\gQEll^\ast_G \otimes_{\Z[q^\pm]} \Z\llpar q \rrpar \xrightarrow{\simeq} \gKU_\fin\llpar q \rrpar^\ast_G\]
    is an isomorphism.
\end{theorem}

This is a highly structured and global version of one of the defining properties of equivariant quasi-elliptic cohomology \cite[(1.1)]{huanquasielliptic}.

\begin{proof}
There is a map $\varphi\colon \gQEll \to \gKU_\fin\llpar q \rrpar$ from \Cref{functorialityoftwoglobal} as the oriented $\P$-divisible group defining $\gKU_\fin\llpar q\rrpar$ is the base-change of $\T_\km$ along $\KU[q^\pm] \to \KU\llpar q \rrpar$. The map (\ref{bcofqell}) above is this map $\varphi$ base-changed to a map of $\E_\infty$-$L \KU\llpar q \rrpar$-algebras. To show this map is an equivalence, it suffices to apply geometric fixed points $\Phi^G$ for each finite abelian group $G$
\begin{equation}\label{bcofqellaftergeofp}
\Phi^G(\gQEll) \otimes_{\KU[q^\pm]} \KU\llpar q \rrpar \simeq \Phi^G\left(\gQEll \otimes_{L\KU[q^\pm]} L\KU\llpar q \rrpar\right) \to \Phi^G(\gKU_\fin \llpar q \rrpar),
\end{equation}
where the first equivalence comes from the fact that $\Phi^G$ is symmetric monoidal and preserves colimits and that there is a natural equivalence $\Phi^G L X \simeq X$ for spectra $X$; see \cite[Pr.4.5.8]{s}. The map (\ref{bcofqellaftergeofp}) is then the base-change equivalence for geometric fixed points of oriented tempered cohomology theories of \cite[Th.5.3.10]{tempered}.\\
    
The ``in particular'' statement essentially follows from a degenerate global Tor-spectral sequence. Due to a lack of references in this direction, we will argue directly. First, consider the morphism of $\E_\infty$-rings $\ku[q^\pm] \to \ku\llpar q \rrpar$ given by the connective cover of $\KU[q^\pm] \to \KU\llpar q \rrpar$. By Lurie's $\E_\infty$-version of Lazard's theorem, see \cite[Th.7.2.2.15]{ha}, we can write $\ku\llpar q \rrpar \simeq \colim f_i$ as a filtered colimit of finite free $\ku[q^\pm]$-modules. Inverting the Bott element $u$ gives us an expression for $\KU\llpar q \rrpar \simeq \colim F_i$ as a filtered colimit of finite free $\KU[q^\pm]$-modules $F_i$. As $L$ commutes with colimits and preserves free modules, we see that $L \KU\llpar q \rrpar$ is a filtered colimit of finite free $L\KU[q^\pm]$-modules
    \[L (\KU\llpar q \rrpar) \simeq L (\colim F_i) \simeq \colim L F_i.\]
    By \cite[Pr.3.4]{firstpreprint}, we see $L F_i$ is a \emph{globally flat} $L \KU[q^\pm]$-module, which allows us to use the collapsed Tor-spectral sequence of \cite[Pr.4.8]{firstpreprint} producing the natural isomorphisms of rings
    \[\pi_\ast^G \left( \gQEll \underset{L\KU[q^\pm]}{\otimes} L F_i\right) \simeq \pi_\ast^G \gQEll \underset{\pi_\ast^e L \KU[q^\pm]}{\otimes} \pi_\ast^e LF_i \simeq \pi_\ast^G \gQEll \underset{\Z[q^\pm, u^\pm]}{\otimes} \pi_\ast F_i.\]
    Combing these facts with the previous paragraph and the fact that filtered colimits commute with $G$-homotopy groups and relative tensor products, we obtain the desired isomorphism of rings
    \[\gQEll^\ast_G \underset{\Z[q^\pm]}{\otimes} \Z\llpar q \rrpar \simeq \pi_{-\ast}^G \gQEll \underset{\pi_{-\ast} \KU[q^\pm]}{\otimes} \pi_{-\ast} \KU\llpar q \rrpar \simeq \colim \left( \pi_{-\ast}^G \gQEll \underset{\pi_{-\ast}\KU[q^\pm]}{\otimes}\pi_{-\ast} F_i \right)\]
    \[\simeq \colim \pi_{-\ast}^G \left( \gQEll \underset{L \KU[q^\pm]}{\otimes} L F_i \right) 
    \simeq \pi_{-\ast}^G \left(\gQEll \underset{L \KU[q^\pm]}{\otimes} L\KU\llpar q \rrpar\right)
    \simeq \gKU_\fin\llpar q \rrpar_G^\ast\]
    for each finite group $G$.
\end{proof}

As $\gQEll$ comes from a tempered cohomology theory, many of the properties of equivariant quasi-elliptic cohomology proven in \cite{huanquasielliptic} can be refined and generalised. We give one more example. For a global space $X$, we write $\gQEll^X$ for the global $\E_\infty$-ring given by $\underline{\Map}(\Sigma^\infty_+X,\gQEll)$ , where $\underline{\Map}$ denotes the internal hom in global spectrum, and $\gQEll^X_e$ for the underlying $\E_\infty$-ring. In particular, $\QEll_e = \KU[q^\pm]$ by construction. 

\begin{cor}\label{kuenneth}
    Let $X$ and $Y$ be $\fin$-global spaces such that $\gQEll_e^X$ is a perfect $\KU[q^\pm]$-module. Then the natural map of $\E_\infty$-$\gQEll_e$-algebras
    \[\gQEll_e^X \otimes_{\gQEll_e} \gQEll_e^Y \to \gQEll_e^{X\times Y}\]
    is an equivalence. Furthermore, if $\gQEll_e^\ast(X)$ is flat as a module over $\gQEll^\ast(\ast)\simeq \Z[q^\pm, u^\pm]$, then the natural map of commutative $\Z[q^\pm, u^\pm]$-algebras
    \[\gQEll^\ast(X) \otimes_{\Z[q^\pm, u^\pm]} \gQEll^\ast(Y) \to \gQEll^\ast(X\times Y)\]
    is an isomorphism. In particular, if $H,G$ are finite groups, the natural map of $\Z[q^\pm,u^\pm]$-algebras
    \[\gQEll^\ast_H \underset{\Z[q^\pm,u^\pm]}{\otimes} \gQEll^\ast_G\xrightarrow{\simeq} \gQEll_{H\times G}^\ast\]
    is an isomorphism
\end{cor}

\begin{remark}
Given a $\fin$-global space $X$, $\gQEll_e^X$ is a perfect $\QEll_e$-module if
\begin{enumerate}
    \item $X$ is equivalent to $\R X'$ for some $\pi$-finite space $X'$ (\cite[Pr.4.7.9]{ec3}), or 
    \item $X$ is equivalent to a global quotient $X''//G$ for a finite group $G$ and a finite $G$-space $X''$ (\cite[Cor.4.7.10]{ec3}).
\end{enumerate}
\end{remark}

\begin{proof}[Proof of \Cref{kuenneth}]
    This follows from \cite[Pr.4.4.7]{ec3} and the ``furthermore'' statement from a degenerate Tor-spectral sequence in $\Sp$. For the ``in particular'' statement, we first note that a consequence of Lurie's \emph{tempered ambidexterity theorem} \cite[Th.7.2.10]{ec3} is that $\gQEll_e^{\B G}$ is perfect as a $\KU[q^\pm]$-module for all finite groups $G$; see \cite[Rmk.10.6]{davidandlennart}. We then compute that the graded $\Z[q^\pm, u^\pm]$-module $\pi_{-\ast}\gQEll_e^{\B G} = \gQEll^\ast_G(\ast)$ to be finite free, see \cite[Lm.3.1]{huanquasielliptic}, which gives us the desired isomorphisms of rings
\[\gQEll^\ast_H \underset{\Z[q^\pm,u^\pm]}{\otimes} \gQEll^\ast_G \simeq \gQEll^\ast(\B H) \underset{\Z[q^\pm,u^\pm]}{\otimes} \gQEll^\ast(\B G) \simeq \gQEll^\ast(\B(H\times G)) \simeq \gQEll_{H\times G}^\ast\]
for finite groups $H,G$.
\end{proof}

\begin{remark}
    One can also define a real version of $\gQEll$ from $\Spec \KU[q^\pm] / C_2$ and show that the natural map to $\gQEll$ is a $C_2$-Galois extension. The follows from \cite{tempered} and will appear in forthcoming work of the first-named author.
\end{remark}

\begin{remark}
As mentioned we may extend quasi-elliptic cohomology to a global cohomology theory on all compact Lie groups by right Kan extension. However, remarks of Rezk suggest that it may be possible to non-trivially extend quasi-elliptic cohomology to a global cohomology theory valued in global spaces with isotropy in compact abelian Lie groups. The $\TT$-fixed points of this theory would be a spectral refinement of the stack denoted $\G_m//q^{\Z}$ in \cite[\textsection VII]{dr} and \cite{rezkquasielliptic}. Such an extension may allow one to more directly approach the derived Tate curve by base-changing this stack to $\KU\llpar q\rrpar$, thus avoiding an appeal to the Serre--Tate theorem.\\

Following the pattern of this article, one may expect that the stack $\G_m//q^{\Z}$ should be given by $\Ell_{\G_m/\KU}(\Cyc(\B \TT))$. Note that we have not even defined the cyclification of $\B\TT$, but we see this as the lesser problem. However, it is not clear how to the authors how to control the result of such a construction. If one could directly construct an elliptic curve in this way, then the uniqueness statement of \Cref{liftoftate} would likely show that this elliptic curve agrees with $\T$.
\end{remark}


\subsection{Two \texorpdfstring{$C_2$}{C2}-equivariant global meromorphic topological \texorpdfstring{$q$}{q}-expansion maps}\label{ssec:real}

Recall the $C_2$-action on $\KU\llpar q \rrpar$ of (\ref{eq:ctwoaction}); by \Cref{complexconjguation}, the nontrivial map associated with this $C_2$-action agrees with the Adams operation $\psi^{-1}$ and is related to Atiyah's \emph{Real $K$-theory} from \cite{ktheoryandreality}. There is also an Adams operations $\psi^{-1}$ acting on $\TMF$, however, this action is trivial; see \cite[Pr.2.15]{heckeontmf}. There are many forms of $\TMF$ with level structure upon which this $C_2$-action is nontrivial. The minimal such example is $\TMF_1(3)$.\\

Recall from \cite[\textsection7]{km}, that $\M_1(3)$ is the moduli stack of elliptic curves $E$ over $\Spec \Z[\tfrac{1}{3}]$ equipped with a point $P$ of exact order $3$. The forgetful map of stacks $p\colon \M_1(3) \to \M_\Ell \times \Spec \Z[\tfrac{1}{3}]$ is finite étale, hence we obtain a stack $\M_1^\ori(3) = (\M_1(3), p^\ast \O_{\M_\Ell^\ori})$ with $\Ga(\M_1^\ori(3)) = \TMF_1(3)$. One can also identify the presheaf on $\E_\infty$-rings corepresented by the stack $\M_1^\ori(3)$ with triples $(E,P,e)$, where $(E,e)$ is an oriented elliptic curve over $R$ and $P\colon \Z/3 \to E(R)$ is a choice of nonzero map in $\Mod_\Z$. We will call $P$ a \emph{point of order $3$}.\\

There is a natural involution $\tau$ on $\M_1^\ori(3)$ defined by the triple $(\e_1(3),-P,-e)$, where $(\e_1(3),e)$ is the pullback of the universal oriented elliptic curve along the forgetful map $\M_1^\ori(3) \to \M_\Ell^\ori$ and $P$ is the universal point of order $3$. This defines a $C_2$-action over $\M_\Ell^\ori$, as the involution $[-1]\colon \e_1(3) \to \e_1(3)$ witnesses the commutativity of the diagram
\[\begin{tikzcd}
    {\M_1^\ori(3)}\ar[rr, "\tau"]\ar[rd, "{(\e_1(3), e)}", swap]    &&  {\M_1^\ori(3)}\ar[dl, "{(\e_1(3),e)}"]  \\
        &   {\M_\Ell^\ori.}  &
\end{tikzcd}\]
One can now apply \Cref{functorialityoftwoglobal} to obtain a $C_2$-action on the global $\E_\infty$-ring $\gTMF_1(3)$, which furthermore is $\gTMF$-linear. The underlying associated Borel $C_2$-spectrum $\TMF_1(3)$ and its cousins $\Tmf_1(3)$ and $\tmf_1(3)$ were extensively analysed in \cite{levelonethree,lennartandhill}.\\

We would now like to show how the derived Tate curve $\T$ can be used to construct $C_2$-equivariant maps of global $\E_\infty$-rings
\[  \ga \colon \gTMF_1(3) \to \gKU\llpar q \rrpar(\omega), \qquad \rho \colon \gTMF_1(3) \to \gKU\llpar q^{1/3} \rrpar[\tfrac{1}{3}],\]
where the codomains will be defined shortly, corresponding to the ramified and unramified cusps of $\overline{\M}_1(3)$.

\subsubsection{The unramified cusp}

Let us write $\Sph(\omega)$ for the finite \'{e}tale extension of $\Sph[\tfrac{1}{3}]$ given by the finite \'{e}tale map of rings $\Z[\tfrac{1}{3}] \to \Z[\tfrac{1}{3},\omega]$ adjoining a primitive $3$\textsuperscript{rd} root of unity; see \cite[Cor.7.5.4.3]{ha} or \cite[Th.3]{svw}. In particular, $\Sph(\omega)$ is a $C_2$-Galois extension of $\Sph[\tfrac{1}{3}]$. Define $\E_\infty$-rings
\[\KU[q^\pm, \omega] = \KU[q^\pm] \otimes \Sph(\omega) \qquad\text{and}\qquad \KU\llpar q \rrpar (\omega) = \KU\llpar q\rrpar \otimes \Sph(\omega),\]
where the latter is also equipped with the diagonal $C_2$-action induced by this Galois action on $\Sph(\omega)$ and the $[-1]$-action of \Cref{eq:ctwoaction}. Let us denote this diagonal involution on $\Spec \KU\llpar q \rrpar(\omega)$ by $\Delta c$.

\begin{theorem}\label{realunramqexpansionmap}
    There is a $C_2$-equivariant map of stacks
    \[\ga\colon \Spec \KU\llpar q \rrpar(\omega) \to \M_1^\ori(3)\]
    inducing a map in $\Fun(BC_2, \CAlg(\Sp_\gl))$ of the form
    \[\ga \colon \gTMF_1(3) \to \gKU\llpar q\rrpar(\omega),\]
    whose $C_2$-fixed points and underlying map fit into a commutative diagram of global $\E_\infty$-rings
    \begin{equation}\label{realunramcommutativediagram}\begin{tikzcd}
        {\gTMF}\ar[d, "{\T}"]\ar[r]  &   {\gTMF_1(3)^{hC_2}}\ar[d, "{\ga^{hC_2}}"]\ar[r]  &    {\gTMF_1(3)}\ar[d, "{\ga}"]\\
        {\gKO\llpar q \rrpar}\ar[r]  &   {\gKU\llpar q \rrpar(\omega)^{h C_2}}\ar[r]  &   {\gKU\llpar q \rrpar(\omega)}
    \end{tikzcd}\end{equation}
    containing the global meromorphic topological $q$-expansion map of \Cref{globalmap} induced by $\T$. We call $\ga$ the \emph{unramified $C_2$-equivariant global meromorphic topological $q$-expansion map}.
\end{theorem}

Note that there is a natural equivalence of $\E_\infty$-rings $\TMF_1(3)^{hC_2} \simeq \TMF_0(3)$ as $\M_1(3) \to \M_0(3)$ is a $C_2$-torsor, so $\TMF_0(3) \to \TMF_1(3)$ is a $C_2$-Galois extension; this follows by étale descent, but also see \cite[Th.7.6]{akhilandlennart}, for example. The spectrum $\KU\llpar q \rrpar(\omega)^{hC_2}$ is a Laurent series version of the spectrum $\KU^\tau$ from \cite{tylerniko} at the end of Remark A.6. One can compute
\[\pi_\ast \KU\llpar q \rrpar(\omega)^{hC_2} \simeq \Z\llpar q \rrpar[\tfrac{1}{3}, (\sqrt{-3}u)^{\pm}].\]

\begin{proof}
    The constant loop map $X \to L X // S^1$ induces a natural transformation from the inclusion of $\Spc_{\ab\fin} \to \Spc_{\gl}$ to $\Cyc$. This induces a map of oriented $\P$-divisible groups
    \begin{equation}\label{inclusionofmultiplicativepdiv} \G_m[\P^\infty] = \mu_{\P^\infty} \to \T_\km \end{equation}
    over $\KU[q^\pm]$, the former being the oriented $\P$-divisible group associated with the base-change of the oriented geometric abelian $\KU$-stack $\G_m$ to $\KU[q^\pm]$. Let us now implicitly invert $3$ everywhere. Over $\Sph(\omega)$, there is a chosen isomorphism $\underline{\Z/3} \simeq \G_m[3] = \mu_3$ given by choosing the element $\omega \in \G_m(\Sph(\omega))$ of order $3$. This gives a choice of point of order $3$ to the finite flat group scheme $\G_m[3]$ over $\Sph(\omega)$. As $\mu_{\P^\infty}$ is the $\P$-divisible group associated with $\G_m$, their $3$-torsion schemes agree; see \cite[Pr.4.3]{elltempcomp}, for example. In particular, we see that any $\E_\infty$-$\Sph(\omega)$-algebra, such as $\KU[q^\pm, \omega]$, comes with a chosen point of order $3$ in $\mu_{\P^\infty}$. The map (\ref{inclusionofmultiplicativepdiv}) then gives $\T_\km$ a chosen point of order 3 over $\KU[q^\pm,\omega]$. Consider the oriented elliptic curve $\T$ now over $\KU\llpar q \rrpar(\omega)$. By \Cref{liftoftate} and the fact that taking associated $\P$-divisible groups commutes with base-change, we see that the $\P$-divisible group associated with $\T$ agrees with the base-change of $\T_\km$ to $\KU\llpar q \rrpar(\omega)$, which equips $\T$ with a choice of point of order $3$ over $\KU\llpar q \rrpar(\omega)$. This provides us with a map of stacks
    \[(\T,\omega,e)\colon \Spec \KU\llpar q \rrpar(\omega) \to \M_1^\ori(3)\]
    over $\M_\Ell^\ori$. We claim this map naturally refines to a $C_2$-equivariant map of stacks over $\M_\Ell^\ori$. This is clear though, as the diagram of stacks over $\M_\Ell^\ori$
    \[\begin{tikzcd}
        {\Spec \KU\llpar q \rrpar(\omega)}\ar[r, "{\Delta c}"]\ar[d, "{(\T,\omega,e)}"] &   {\Spec \KU\llpar q \rrpar(\omega)}\ar[d, "{(\T,\omega,e)}"] \\
        {\M_1^\ori(3)}\ar[r, "\tau"]    &   {\M_1^\ori(3)}
    \end{tikzcd}\]
    commutes: both composites are naturally equivalent to the triple $(\T,-\omega,-e)$ over $\M_1^\ori(3)$, as $c\Delta$ sends $\omega$ to $-\omega$ from the Galois action on $\Sph(\omega)$ and $e$ to $-e$ from the complex conjugation action on $\KU\llpar q \rrpar$. It is important here to use the \emph{diagonal} $C_2$-action on $\KU\llpar q \rrpar(\omega)$ else the upper composite above would be either the triple $(\T,\omega,-e)$ or $(\T,-\omega, e)$. One then obtains the $C_2$-equivariant morphism of global $\E_\infty$-rings by applying \Cref{functorialityoftwoglobal} to the $C_2$-equivariant map of stacks above; here $\x=\M_\Ell^\ori$. The diagram (\ref{realunramcommutativediagram}) then follows by taking $C_2$-fixed points and underlying global $\E_\infty$-rings.
\end{proof}

\begin{remark}[Commutativity with Adams operations]\label{commutativitywithadamsoperations}
    Recall the Adams operations $\psi^k$ on $\KU\llpar q \rrpar[\tfrac{1}{k}]$ from \Cref{rmk:adamsoperations}. There are also Adams operations $\psi^k$ on $\TMF_1(3)[\tfrac{1}{k}]$ from \cite[Df.2.1]{heckeontmf}. It is clear from the descriptions of both of these Adams operations that they are $C_2$-equivariant; this is discussed in \Cref{rmk:adamsoperations} and \cite[Th.C]{heckeontmf}, respectively. We also claim that the map of $\E_\infty$-rings
    \[\ga \colon \TMF_1(3)[\tfrac{1}{k}] \to \KU\llpar q\rrpar(\omega)[\tfrac{1}{k}]\]
    commutes with the Adams operation $\psi^k$ on both sides. To see this, just as in the proof of \Cref{realunramqexpansionmap}, it suffices to show that $\psi^k$ on $\TMF_1(3)$ is defined in \cite[Th.C]{heckeontmf} is equivalent to the global sections of the map $\M_1^\ori(3) \to \M_1^\ori(3)$ defined by the triple $(\e_1(3),[k](P), [k](e))$ where $[k]\colon \e_1(3) \to \e_1(3)$ is the $k$-fold multiplication map. Notice this does not define a $\TMF$-linear map as $[k]$ is not an \emph{isomorphism} of elliptic curves for $k\neq \pm1$. 
\end{remark}

\subsubsection{The ramified cusp}
To construct our second meromorphic topological $q$-expansion map, we start with three lemmata that better describe $\T_\km$ as a $\P$-divisible group over $\KU[q^\pm]$.

\begin{lemma}\label{lm:kmses}
    The natural map $\mu_{\P^\infty} \to \T_\km$ of preoriented $\P$-divisible groups over $\KU[q^\pm]$ induced by the inclusion of constant loops $X \to LX // S^1$ has cokernel $\underline{\Q/\Z}$. In particular, for each $N\geq 2$ we have a short exact sequence of finite flat commutative group schemes over $\KU[q^\pm]$
    \begin{equation}\label{eq:kmses}
        0 \to \mu_N \to \T_\km[N] \to \underline{\Z/N} \to 0.
    \end{equation}
\end{lemma}

\begin{proof}
    First, note that the natural inclusion $\mu_{\P^\infty} \to \T_\km$ is a monomorphism of $\P$-divisible groups \`{a} la \cite[Not.2.7.11]{ec3}. Indeed, over each $p$-local component $\mu_{p^\infty}\to (\T_\km)_{(p)}$ this induces a monomorphism of abelian groups when evaluated on discrete $\tau_{\geq 0}\KU[q^\pm]=\ku[q^\pm]$-algebras as this is true classically. By \cite[Pr.2.4.8]{ec2}, we see that the category of $\P$-divisible groups over $\KU[q^\pm]$ admits a quotient $\T_\km/\mu_{\P^\infty}$. By \cite[Rmk.2.5.4]{ec2}, to see that this quotient is \'{e}tale it suffices to check this on all residue fields of $\pi_0\KU[q^\pm]\simeq \Z[q^\pm]$ and this is again clear as it is true classically; see \cite[(8.7.2.1)]{km}. As the category of \'{e}tale $\P$-divisible groups over $\KU[q^\pm]$ is equivalent to the category of \'{e}tale $\P$-divisible groups over $\Z[q^\pm]$ by taking $\pi_0$, we can also identify this quotient as $\underline{\Q/\Z}$, as desired; see \cite[Pr.2.5.9]{ec2}.
\end{proof}

\begin{lemma}\label{lm:kmsplitting}
    There is a faithfully flat $\KU[q^\pm]$-algebra $\SS$ and a map of $\P$-divisible groups $\underline{\Q/\Z}\to \T_\km$ over $\SS$ inducing a splitting
    \[\mu_{\P^\infty} \oplus \underline{\Q/\Z} \simeq \T_\km\]
    of $\P$-divisible groups over $\SS$. In fact, $\SS$ is the filtered colimit of finite faithfully flat $\KU[q^\pm]$-algebras $\SS_N$ together with maps of finite flat groups schemes $\underline{\Z/N} \to \T_\km[N]$ over $\SS_N$ inducing a splitting
    \[\mu_N \oplus \underline{\Z/N} \simeq \T_\km[N]\]
    of finite flat group schemes over $\SS_N$.
\end{lemma}

If we work over $\SS_3$, then we see that $\T_\km$ has another interesting $3$-torsion point given by the above splitting.

\begin{proof}
    The first statement is simply an application of \cite[2.7.15]{ec3}, and the statement at each finite level was achieved throughout the course of its proof. Indeed, the finite faithfully flat $\KU[q^\pm]$-algebra $\SS_N$ corepresents the functor $X^\circ_N\colon \CAlg_{\KU[q^\pm]} \to \Spc$ given by the subfunctor of
    \[X_N(\A) = \Map_{\Mod_\Z}(\Z/N, \T_\km(\A))\]
    spanned by those maps which induces equivalences $\underline{\Z/N}\simeq (\T_\km/\mu_{\P^\infty})[N]$ of finite flat group schemes over $\A$.
\end{proof} 

\begin{lemma}
    The finite faithfully flat $\E_\infty$-$\KU[q^\pm]$-algebra $\SS_N$ of \Cref{lm:kmsplitting} has homotopy groups
    \[\pi_\ast \SS_N\simeq (\pi_0 \SS_N)[u^\pm] \simeq \Z[q^{\pm 1/N}][u^\pm],\qquad |u|=2.\]
\end{lemma}

\begin{proof}
    By \cite[Rmk.2.7.17]{ec3}, we see that $\pi_0 \SS$ is the universal splitting algebra for $\T_\km^\heartsuit$ over $\Z[q^\pm]$ which is alluded to in \cite[(8.7.4.1)]{km} but not given an explicit form. The same reasoning shows that $\pi_0 \SS_N$ is the universal splitting algebra for $\T_\km^\heartsuit[N]$ over $\Z[q^\pm]$, which is given the explicit form $\Z[q^{\pm 1/N}]$ in \cite[(8.7.1.6)]{km}, which gives our $\pi_0$-computation. For the $\pi_\ast$-computation, we use the fact that $\SS_N$ is flat over $\KU[q^\pm]$.
\end{proof}

We now need to give $\SS_N$ a $C_2$-action. This is a little more complicated than the diagonal $C_2$-action on $\KU\llpar q \rrpar(\omega)$ as that is simply a coproduct in $\Fun(BC_2,\CAlg)$.

\begin{construction}
    Recall that $\KU[q^{\pm}]$ is canonically a $C_2$-ring, whose action we denote by $c$, and note that this induces a $C_2$-action on the category $\CAlg_{\KU[q^{\pm}]}$. Now as previously observed, $\T_{\km}(\KU[q^{\pm}])$ is naturally equivalent to $\T_{\km}(c^* \KU[q^{\pm}]) = (c^\ast \T_\km)(\KU[q^{\pm}])$. In turn this induces a natural equivalence $X_N^\circ(-)\simeq X_N^\circ(c^*(-))$. In fact the identification above is $C_2$-equivariant, meaning that $\SS_N$ is canonically an object of $\CAlg_{\KU[q^{\pm}]}^{h C_2}$. Note now that $\SS_N$ inherits a second action by $(\Z/N)^\times$-action via precomposition on $X^\circ_N$. Moreover, this clearly commutes with the $C_2$-action constructed before by naturality. In particular, we obtain $\SS_3$ as an object of $\Fun(BC_2, \CAlg_{\KU[q^{\pm}]}^{hC_2})$. Now recall that the category of $\E_\infty$-algebras under $\KU[q^{\pm}]$ is a slice category. Since passing to slices commutes with limits we find that in general, given a $G$-equivariant object $X$ in a category $\mathcal{C}$,
    \[(\mathcal{C}_{X/})^{hG} \simeq (\mathcal{C}^{hG})_{X/} \simeq \Fun(BG,\mathcal{C})_{X/},\] 
    where the final equivalence is justified by the fact that $\mathcal{C}$ itself has a trivial $G$-action. In particular, we obtain $\SS_3$ as an object of $\Fun(B C_2, \Fun(B C_2, \CAlg)_{\KU[q^{\pm}]/})$. Applying the previous observation in the other direction gives $\SS_3$ as an object of $\Fun(B C_2\times B C_2, \CAlg)_{\KU[q^{\pm}]/}$, where now $\KU[q^{\pm}]$ is a $C_2\times C_2$-object via complex conjugation in the second factor and the identity in the first. Passing to the diagonal copy of $C_2$, we obtain a map $\KU[q^{\pm}]\to \SS_3$ which is $C_2$-equivariant, where the source is given the complex conjugation action. Let us also write $c\Delta$ for this $C_2$-action on $\SS_3$.
\end{construction}

\begin{mydef}
    Let $\KU\llpar q^{1/3}\rrpar$ be the $\E_\infty$-$\KU\llpar q \rrpar$-algebra in $\CAlg^{BC_2}$ defined by
    \[\KU\llpar q^{1/3}\rrpar = \KU\llpar q \rrpar \underset{\KU[q^\pm]}{\otimes} \SS_3.\]
\end{mydef}

Just as it was important to give $\KU\llpar q \rrpar(\omega)$ the diagonal $C_2$-action to prove \Cref{realunramqexpansionmap}, the above $C_2$-action on $\KU\llpar q^{1/3}\rrpar$ is crucial in proving the following.

\begin{theorem}\label{realramfqexpansionmap}
    There is a $C_2$-equivariant map of stacks
    \[\rho\colon \Spec \KU\llpar q^{1/3}\rrpar \to \M_1^\ori(3)\]
    inducing a map in $\Fun(BC_2, \CAlg(\Sp_\gl))$ of the form
    \[\rho \colon \gTMF_1(3) \to \gKU\llpar q^{1/3} \rrpar[\tfrac{1}{3}],\]
    whose $C_2$-fixed points and underlying map fit into a commutative diagram of global $\E_\infty$-rings
    \begin{equation}\label{realramcommutativediagram}\begin{tikzcd}
        {\gTMF}\ar[d, "{\T}"]\ar[r]  &   {\gTMF_1(3)^{hC_2}}\ar[d, "{\rho^{hC_2}}"]\ar[r]  &    {\gTMF_1(3)}\ar[d, "{\rho}"]\\
        {\gKO\llpar q \rrpar}\ar[r]  &   {\gKU\llpar q^{1/3} \rrpar[\tfrac{1}{3}]^{hC_2}}\ar[r]  &   {\gKU\llpar q^{1/3} \rrpar[\tfrac{1}{3}]}
    \end{tikzcd}\end{equation}
    containing the global meromorphic topological $q$-expansion map of \Cref{globalmap} induced by $\T$. We call $\rho$ the \emph{ramified $C_2$-equivariant global meromorphic topological $q$-expansion map}.
\end{theorem}

In this case, notice that $(\KU\llpar q^{1/3}\rrpar[\frac{1}{3}])^{hC_2}$ can be computed to be $\KO\llpar q^{1/3}\rrpar[\tfrac{1}{3}]$ using the usual homotopy fixed point spectral sequence for $\KO \simeq \KU^{hC_2}$; we can even import the key $d_3$-differential along the $C_2$-equivariant map of $\E_\infty$-rings $\KU \to \KU\llpar q^{1/3}\rrpar[\tfrac{1}{3}]$.

\begin{proof}
    The splitting of \Cref{lm:kmsplitting} induces a chosen inclusion $s\colon \underline{\Z/3} \to \T_\km[3]$ over $\KU[q^{\pm 1/3}]$. As in the proof of \Cref{realunramqexpansionmap}, as the $3$-torsion of $\T$ over $\KU\llpar q^{1/3}\rrpar$ agrees with the $3$-torsion of the base-change of $\T_\km$, we obtain a chosen point of order $3$ given by the image of the generator $1$ of $\Z/3\Z$ inside $\T[3] \subseteq \T$ over $\KU\llpar q^{1/3}\rrpar$. This produces the desired map of stacks
    \[(\T,s,e) \colon \Spec \KU\llpar q^{1/3}\rrpar[\tfrac{1}{3}] \to \M_1^\ori(3)\]
    over $\M_\Ell^\ori$. The only reason we need to invert $3$ here is to work over $\M_1^\ori(3)$ where $3$ is already inverted. Again, as in the proof of \Cref{realunramqexpansionmap}, this map is naturally $C_2$-equivariant. Indeed, this comes from the fact that the diagram of $\E_\infty$-rings
    \[\begin{tikzcd}
        {\Spec \KU\llpar q^{1/3} \rrpar[\tfrac{1}{3}]}\ar[r, "{\Delta c}"]\ar[d, "{(\T,s,e)}"] &   {\Spec \KU\llpar q^{1/3} \rrpar[\tfrac{1}{3}]}\ar[d, "{(\T,s,e)}"] \\
        {\M_1^\ori(3)}\ar[r, "\tau"]    &   {\M_1^\ori(3)}
    \end{tikzcd}\]
    naturally commutes, as both composites are represented by the triple $(\T,-s,-e)$. Again, the $(\Z/3)^\times$-action on $\SS_3$ sends $s$ to $-s$ and the $C_2$-action on $\KU\llpar q \rrpar$ sends $e$ to $-e$. This $C_2$-equivariant map of stacks induces the desired $C_2$-equivariant map of global $\E_\infty$-rings by \Cref{functorialityoftwoglobal}. The desired commutative diagram of $\E_\infty$-rings (\ref{realramcommutativediagram}) then follows by taking underlying global $\E_\infty$-rings and $C_2$-fixed points.
\end{proof}

\begin{remark}[Connection to geometric fixed points]
    It feels more natural to us to replace $\TMF_1(3)$ with the geometric fixed points $\Phi^{\Z/3}\gTMF_\fin$, or more generally, with the fin-global $\E_\infty$-ring $\gTMF_\fin^{\Phi \Z/N}$. These geometric fixed points are an integral refinement of $\TMF_1(3)$, hence provide a more powerful map without needing to invert $3$ over the ramified cusp. Here, we use \cite[Th.F]{tempered} that shows that $\Phi^{\Z/3} \gTMF_\fin$ is the global sections of a substack $\underline{\mathrm{Inj}}(\underline{\Z/N}^\vee, \e[\P^\infty])$ of $\HOM(\underline{\Z/N}^\vee, \e[\P^\infty])$ over $\M_\Ell^\ori$ spanned by those ``injective'' homomorphisms. This then defines a fin-global $\E_\infty$-ring $\gTMF_\fin^{\Phi \Z/3}$ whose underlying $\E_\infty$-ring is $\Phi^{\Z/3} \gTMF_\fin$. The splitting of $\T_\km[N]$ over $\SS_N$ then gives the desired map of stacks
    \[\Spec \SS_N \to \underline{\mathrm{Inj}}(\underline{\Z/N}^\vee, \T_\km)\]
    over $\KU[q^\pm]$, and hence a map of fin-global $\E_\infty$-rings $\gQEll^{\Phi \Z/3} \to \SS_N^\gl$. By base-change and precomposing with the map of geometric fixed points induced by $\gTMF_\fin \to \gQEll$, we obtain the desired map of fin-global $\E_\infty$-rings:
    \[\gTMF_\fin^{\Phi \Z/N} \to \gKU_\fin\llpar q \rrpar^{\Phi \Z/N} \simeq \gQEll^{\Phi \Z/N} \underset{\KU[q^\pm]}{\otimes} \KU\llpar q \rrpar \to  \SS_N^\gl \underset{\KU[q^\pm]}{\otimes} \KU\llpar q \rrpar = \KU_\fin\llpar q^{1/N}\rrpar.\]
\end{remark}

\begin{remark}[$C_2$-equivariant holomorphic topological $q$-expansion maps]
    Using these $C_2$-equivariant maps $\ga$ and $\rho$ of \Cref{realunramqexpansionmap,realramfqexpansionmap}, one could define a $C_2$-action on $\Tmf_1(3)$ using the pullback in $C_2$-equivariant $\E_\infty$-rings
    \[\begin{tikzcd}
        {\Tmf_1(3)}\ar[r]\ar[d, "{(\overline{\rho}, \overline{\ga})}"] &   {\TMF_1(3)}\ar[d, "{(\rho,\ga)}"]    \\  
        {\KU\llbracket q^{1/3} \rrbracket[\frac{1}{3}] \times \KU\llbracket q \rrbracket(\omega)}\ar[r]   &   {\KU\llpar q^{1/3} \rrpar[\frac{1}{3}] \times \KU\llpar q \rrpar(\omega).}
    \end{tikzcd}\]
    It is likely that this $C_2$-action on $\Tmf_1(3)$ matches that of \cite{lennartandhill}. If this is the case, then such a pullback would provide two $C_2$-equivariant holomorphic topological $q$-expansion maps $\overline{\rho}$ and $\overline{\ga}$. Such $q$-expansion maps could be incredibly useful in computations.
\end{remark}

\addcontentsline{toc}{section}{References}
\scriptsize

\bibliography{references}
\bibliographystyle{alpha}

\end{document}